\documentclass[11pt]{amsart}

\usepackage{amssymb, amsmath, amscd, amsthm, color, epsfig}

\usepackage[all]{xy}          
\xyoption{dvips}              

\newcommand{\R}{{\mathbb R}}
\newcommand{\Z}{{\mathbb Z}}

\newcommand{\deldot}{{\Delta^{\bullet}}}
\newcommand{\Xdot}{{X^{\bullet}}}

\newcommand{\cd}[1]{\mathcal{D}^{c}_{#1}}
\newcommand{\Tstar}[1]{T^{*}_{#1}}

\newcommand{\K}{{\mathcal{K}}}
\newcommand{\W}{{\mathcal{W}}}
\newcommand{\V}{{\mathcal{V}}}

\newcommand{\FM}[2]{F[#1, #2]}

\newcommand{\HO}{{\mathcal{H}}}

\theoremstyle{plain}
\newtheorem{thm}{Theorem}[section]
\newtheorem{prop}[thm]{Proposition}
\newtheorem{lemma}[thm]{Lemma}
\newtheorem{cor}[thm]{Corollary}

\theoremstyle{definition}
\newtheorem{definition}[thm]{Definition}

\theoremstyle{remark}
\newtheorem*{rem}{Remark}
\newtheorem*{rems}{Remarks}

\newcommand{\refT}[1]{Theorem~\ref{T:#1}}
\newcommand{\refC}[1]{Corollary~\ref{C:#1}}
\newcommand{\refP}[1]{Proposition~\ref{P:#1}}

\begin{document}

\title{Finite type knot invariants and calculus of functors}
\author{Ismar Voli\'c}
\address{Department of Mathematics, University of Virginia,
Charlottesville, VA}
\email{ismar@virginia.edu}
\urladdr{http://www.people.virginia.edu/\~{}iv2n}

\begin{abstract}
We associate a Taylor tower supplied by calculus of the embedding
functor to the space of long knots and study its cohomology spectral
sequence.  The combinatorics of the spectral sequence along the line
of total degree zero leads to chord
diagrams with relations as in finite type knot theory.
We show that the spectral sequence collapses along this line and
that
the Taylor tower represents a universal
finite type knot invariant.
\end{abstract}

\maketitle


\section{Introduction}

\setcounter{section}{1}

In this paper we bring together two fields, finite type knot
theory and Goodwillie's calculus of functors, both of which have
received considerable attention during the past ten years.
Goodwillie and his collaborators have known for some time, through
their study of the embedding functor, that there should be a
strong connection between the two \cite{GKW, GW}.  This is indeed
the case, and we will show that a certain tower of spaces in the
calculus of embeddings serves as a classifying object for finite
type knot invariants.

Let $\K$ be the space of \emph{framed long knots}, i.e. embeddings of
$\R$ into $\R^{3}$ which agree with some fixed linear inclusion of
$\R$ into $\R^{3}$ outside a compact set, and come equipped with a
choice of a nowhere-vanishing section of the normal bundle.
Alternatively, we may think of $\K$ as the space of framed \emph{based} knots
in $S^{3}$, namely framed maps of the unit interval $I$ to $S^{3}$ which
are embeddings except at the endpoints.  The endpoints are mapped to
some basepoint  in $S^{3}$ with the same derivative. It is not hard to
see that these two versions of $\K$ are homotopy equivalent.

Now consider (framed) \emph{singular} long knots which are
(framed) long knots with a finite number of transverse double
points. We can extend any knot invariant $V$ to such knots by
applying it to all knots obtained by resolving the double points
of a singular knot in two possible ways, with appropriate signs.
If $V$ vanishes on knots with more than $n$ double points, then it
is of \emph{type $n$}.

As will be explained in more detail in
\S\ref{S:FinTypeInv's}, it turns out that $\V_{n}$, the set of all type $n$ invariants,
 is closely related to the space generated by
chord diagrams with $n$ chords modulo a certain relation.  In fact,
if $\W_{n}$ denotes the dual of this space, then Kontsevich
\cite{Kont} proves
\begin{thm}\label{T:Kontsevich}
$\V_{n}/\V_{n-1}\cong\W_{n}.$
\end{thm}
\noindent
An alternative proof of this theorem uses Bott-Taubes configuration
space integrals which play a crucial role in the proofs of our
results.  These are briefly described in \S\ref{S:B-TIntegrals}, but
more details can be found in \cite{BT, Vo2}.

The other ingredient we need is the \emph{Taylor
tower for $\K$} defined in some detail in \S\ref{S:TaylorTower}.
Following work of Weiss \cite{We},
one can construct spaces $\HO_{r}$
for each $r>1$,
defined as homotopy limits of certain diagrams of ``punctured knots,''
namely embeddings of $I$ in $S^{3}$ as before but with some number of subintervals
removed. There are canonical maps $\K\to\HO_{r}$ and
$\HO_{r}\to\HO_{r-1}$, and these combine to yield the Taylor tower.

To study the tower, we introduce in \S\ref{S:CosimplicialModel} a
cosimplicial space $\Xdot$ whose $r$th partial totalization is
precisely the stage $\HO_{r}$ (\refT{Tot=holim}).  This
equivalence was shown by Sinha \cite{Dev}.  The advantage of the
cosimplicial point of view is that $\Xdot$ comes equipped with a
(second quadrant) cohomology spectral sequence. This spectral
sequence computes the cohomology of the total complex $T_{r}^{*}$
of the double complex obtained by applying cochains to the
truncated cosimplicial space.  Alternatively, one can pass to real
linear combinations generated by the spaces in the cosimplicial
space and take the $r$th partial totalization to get $T_{r}^{*}$
again.  It is the latter description of this algebraic analog of
the Taylor tower that we use for our proofs
 (see \S\ref{S:TotMapsToFinType}
for exact description of the stages $T_{r}^{*}$), because sums of
resolutions of singular knots can be mapped into such a tower. We
can now state the main result:
\begin{thm}\label{T:IntroMainTheorem1}
$H^{0}(T^{*}_{2n})\cong \V_{n}.$
\end{thm}
\noindent
Moreover, all finite type invariants factor through the algebraic Taylor
tower for
$\K$, as the above isomorphism is induced by a map which is essentially a
collection of evaluation maps (see \refP{evMapstoTot}).
Unfortunately, it is not known that the algebraic stages and
ordinary stages of the Taylor tower have the same homology,
although this is believed to be the case (Sinha \cite{Dev} has shown this
for any space of embeddings of one manifold in another,
as long as the codimension is at least 3).

The starting point in the proof of the above theorem is \refP{CD},
which states that the groups on the diagonal of the $E_{1}$ term of the spectral sequence
can be associated to chord diagrams.  We then show
\begin{itemize}
\item $E^{-2n, 2n}_{\infty}\cong
H^{0}(T^{*}_{2n})/H^{0}(T^{*}_{2n-1})$\ \ \ (\S\ref{S:Iso's&Surj'sInTower}),
\item $E^{-2n, 2n}_{2}\cong \W_{n}$\ \ \ (\refP{E2}),
\item There is a commutative diagram
$$
\xymatrix{ E^{-2n, 2n}_{\infty} \ar@{^{(}->}[r] \ar@{^{(}->}[d] &
\V_{n}/\V_{n-1} \ar@{^{(}->}[d]  \\
E^{-2n, 2n}_{2} \ar[r]^{\cong} & \W_{n}}
$$
(\refP{TotMapsToFinType}).
\end{itemize}
The proof of the first statement is a direct
computation of the first differential in the spectral sequence.
  The other parts follow from
considering how a singular knot maps into the algebraic Taylor tower.

We then use \refT{B-TMainTheorem} to finish the proofs in
\S\ref{S:Universality}.  This theorem essentially states that Bott-Taubes
integration fits as a diagonal map from $\W_{n}$ to $E^{-2n,
2n}_{\infty}$ in the above square and makes the top resulting
triangle commute.
It follows that all the maps in the diagram are isomorphisms over the
rationals (Bott-Taubes integrals are the only reason we cannot
state our results over the integers).  In
addition to deducing \refT{IntroMainTheorem1}, we also have an
alternative proof of \refT{Kontsevich}, as well as
\begin{thm}\label{T:Collapse}
The spectral sequence collapses on the diagonal at $E_{2}$.
\end{thm}
Thus calculus of the embedding functor provides a new point of view
on finite type knot theory.
Some further questions as well as a brief discussion of the potential
importance of this point
of view can be found in \S\ref{S:Future}.

\subsection{Acknowledgements}  I am deeply endebted to my thesis
advisor, Tom Goodwillie, whose knowledge is only matched by his desire
to share it.  Dev Sinha's help over the years has also been invaluable.
I am thankful to Pascal Lambrechts and Greg Arone for many
comments and suggestions.

    \section{Finite type knot invariants}\label{S:FinTypeInv's}

Here we recall some basic features of finite type knot theory.
  More details can be
found in \cite{BN, BN2}.

\vskip 5pt
\noindent
As mentioned earlier, a \emph{singular long knot}
is a long knot except for a finite number
of double points.  The tangent vectors at the double points are required
to be independent.  A knot with $n$ such
self-intersections is called \emph{$n$-singular}.

One can now use the
 \emph{Vassiliev skein relation}, pictured in
Figure \ref{F:VassilievSkein}, to extend any knot invariant
$V$ to singular knots.

\begin{figure}[h]
\begin{center}
\input{skeinrelation.pstex_t}
\caption{Vassiliev skein relation}\label{F:VassilievSkein}
\end{center}
\end{figure}
\noindent
The drawings mean that the
three knots only differ locally in one crossing.  A $n$-singular knot
thus produces
$2^{n}$ resolutions.  The sign convention ensures that the order in
which we resolve the singularities does not matter.

\begin{definition}
\label{D:FiniteTypeInvariant}
$V$ is a \emph{(finite, or Vassiliev) type $n$ invariant} if it
vanishes identically on singular knots with $n+1$ self-intersections.
\end{definition}

Let $\V$ be the collection of all finite type invariants and let
$\V_{n}$ be the set of type $n$ invariants.  It is easy to see, for example,
that $\V_{0}$ and $\V_{1}$ (for unframed knots) both contain only the constant functions
on $\K$.  Also immediate is that
$\V_{n}$ contains $\V_{n-1}$.

Another quick consequence of the definition is that the value of a
type $n$ invariant on an $n$-singular knot only
depends on the placement of its singularities.  This is because if
two $n$-singular knots differ only in the embedding, the
difference of $V\in \V_{n}$ evaluated on one and the other is the
value of $V$ on some $(n+1)$-singular knots (since one can get from one
knot to the other by a sequence of crossing changes).  But $V$ by
definition vanishes on such knots.

It thus follows that the value of $V$ on $n$-singular knots is
closely related to the following objects:

\begin{definition} A \emph{chord diagram of degree $n$} is an oriented interval with
$2n$ paired-off points on it, regarded up to orientation-preserving
diffeomorphisms of the interval.
\end{definition}
The pairs of points can be thought of
as prescriptions for where the singularities on the knot should
occur, while how the rest of the interval is embedded is immaterial.

Let $CD_{n}$ be the set of all chord diagrams with $n$ chords.
If $D$ is an element of $CD_{n}$, and if
$K_{D}$ is any $n$-singular knot with singularities as prescribed by
$D$, by observations above we have a map
\begin{equation}\label{E:EasyMap}
\V_{n}\longrightarrow \{f\colon \R[CD_{n}]\to\R\}
\end{equation}
given by
\begin{equation}\label{E:EasyDirection}
f(D)=V(K_{D})
\end{equation}
and extending linearly.
The kernel is by definition $\V_{n-1}$.

\begin{figure}[h]
\begin{center}
\input{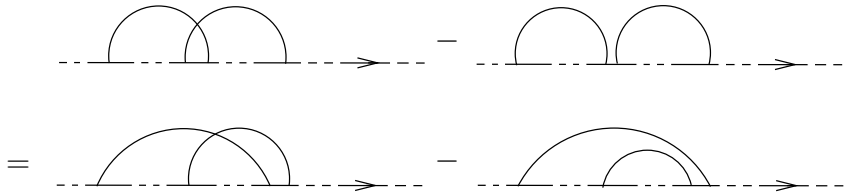}
\caption{$4T$ relation}\label{F:4T}
\end{center}
\end{figure}

Now let $\cd{n}=\R[CD_{n}]/4T$, where $4T$ (or \emph{four-term})
relation is given in  Figure \ref{F:4T}.
The diagrams differ
only in chords indicated; there may be more
chords with their endpoints on the dotted segments, but they are the same for all four
diagrams.  This is a natural relation to impose since it simply
reflects the fact that moving a strand of an $(n-1)$-singular
knot around one of the singularities gets us back to the original
position while introducing four
$n$-singular knots on its circular route.

\begin{rem}
In most literature on finite type theory, one more
relation besides the $4T$ is imposed.  This $1T$ (\emph{one-term)
relation} sets the value of any weight system on a chord diagram with
an \emph{isolated chord}, i.e. a chord not intersected by any other
chords, to be zero.
However, since we will be considering framed knots,
the $1T$ relation cannot be imposed.  This is because
the two resolutions of a singularity coming from an isolated chord
are \emph{not framed isotopic}.
The consequence of having to consider framed knots is simply that the
number of finite type invariants is somewhat larger.  One now
 gets a genuine type 1 invariant,
the framing number.  The $n$th power of the framing
number as well as its product with any type $n-1$ invariant give
additional type $n$ invariants.
\end{rem}

Let $\W_{n}$, the space of \emph{weight systems of degree $n$}, be
the dual of $\cd{n}$. It turns out that $\W_{n}$ is all there is
to the image of the map in \eqref{E:EasyMap}.  As stated in
Introduction, we have the following important theorem due to
Kontsevich \cite{Kont}:

\begin{thm}
$\V_{n}/\V_{n-1}\cong\W_{n}.$
\end{thm}




  The inverse of the map in (\ref{E:EasyMap}) is given by the
famous \emph{Kontsevich Integral}, so it is an example of a \emph{universal
finite type invariant}.  Alternatively, one can use
Bott-Taubes configuration space integrals which we describe in the next
section.

    \section{Bott-Taubes integrals}\label{S:B-TIntegrals}

What follows is a brief outline of how configuration space integrals
can be used to produce a universal finite type invariant.  The
integration techniques mentioned here first arose in Chern-Simons
perturbation theory (see \cite{Lab} for the history of the theory),
but Bott and Taubes \cite{BT} were first to present and develop the
configuration space integrals
in a physics-free way.

Generalizing the computation of the linking number, they
consider configurations of $2n$ points on a knot and use a labeled chord
diagram $D$ as a prescription for constructing $n$ maps to the
product of 2-spheres.

More precisely, let $F(k,N)$ denote the configuration space of $k$
ordered points in a
manifold $N$:
$$
F(k,N)=\{(x_{1},\ldots,x_{k})\in N^{k}, \ \ x_{i}\neq x_{j}.
\text{ for } i\neq j\}.
$$
Let $K$ be a knot.  Given a configuration in
$F(2n,\R)$, $K$ can be evaluated on the $2n$ points to yield
a configuration in $F(2n,\R^3)$.  We denote this map by $ev_K$.
Now let $p_{i}$ and $p_{j}$ be two points in the image of $ev_K$
whose counterparts, $i$ and $j$, are two vertices connected by a
chord in $D$.
The composition of interest is then
\begin{equation}\label{E:IntroComposition}
\xymatrix{
F(2n,\R)\times \K \ar[r]^(0.55){ev_K} & F(2n,\R^{3})
\ar[rr]^(0.55){\frac{p_{j}-p_{i}}{|p_{j}-p_{i}|}} & &
S^{2}.
}
\end{equation}
There are as many maps to $S^2$ as there are chords in $D$.  Their
product can be used for pulling back $n$ standard unit
volume forms from the product of spheres to $F(2n,\R)\times\K$.

If one started with two disjoint knots and a configuration space of
one point on
each of them, the pushforward of the resulting form to the
space of disjoint embeddings of two circles would then yield an
invariant of 2-component links, the linking number.
However, the situation for one knot is not as simple.  The
configuration space $F(2n,\R)$ is not compact, so that the
pushforward to $\K$ may not converge.  What
is required, as it turns out, is a useful compactification of the
configuration spaces appearing in \eqref{E:IntroComposition}.
\vskip 5pt
\noindent
The first construction of the correct compactification is due to Fulton and
MacPherson \cite{FM} (Bott and Taubes use a modification by Axelrod and Singer
\cite{AS}).  At the heart of this compactification is the
blowup along each diagonal of the ordered product of $k$ copies of $N$.
However, we give here an alternative definition due to Sinha \cite{Dev1} which
does not involve blowups and is thus perhaps more accessible.

Assuming $N$ is embedded in a Euclidean space of
dimension $m$, let $i$ be the inclusion of $F(k,N)$ in
$N^{k}$, $\pi_{ij}$ the map to $S^{m-1}$ given by the
normalized difference of points $p_{j}$ and $p_{i}$, and $s_{ijk}$
the map to $[0,\infty]$ given by $|p_{i}-p_{j}|/|p_{i}-p_{k}|$.
\begin{definition}\label{D:DevCompactification}
Let $\FM{k}{N}$ be the closure of the image of $F(k,N)$ in
$N^{k}\times (S^{m-1})^{k\choose 2}  \times [0,\infty]^{k\choose 3}$
under the map $i\times \pi_{ij} \times s_{ijk}$.
\end{definition}

This compactification is a stratified manifold
(manifold with corners), whose stratification is determined by the
rates at which configuration points are approaching each other.  In
particular, a point in a codimension one stratum is determined by
some
number of configuration points colliding at the same time.

The most useful feature of this construction is that the
directions of approach of the colliding points are kept track of.
This allows Bott and Taubes to rewrite \eqref{E:IntroComposition} as
\begin{equation}\label{E:IntroComposition2}
\FM{2n}{\R}\times\K\longrightarrow \FM{2n}{\R^{3}}\longrightarrow S^{2}
\end{equation}
and they show that the product of the 2-forms which are pulled back
from the
spheres extends
smoothly to the boundary of $\FM{2n}{\R}$ (for details, also see
\cite{Vo2}).  The advantage is that now one can
produce a function on the space of knots $\K$ by integrating the
resulting $2n$-form along the compact fiber $\FM{2n}{\R}$ of the
projection
$$\pi\colon \FM{2n}{\R}\times \K\longrightarrow\K.
$$
More precisely, let $\omega$ be the product of the volume forms on
$(S^{2})^{n}$ and let $h_{D}$ be the product of the compositions in
\eqref{E:IntroComposition2}.  What has just been described
 is a function on $\K$, which we denote by $I(D,K)$, given
by
\begin{equation}\label{E:IntroPushforward}
I(D,K)=\pi_{*}(h_{D}^{*}\omega).
\end{equation}
The question now is whether $I(D,K)$ is a closed 0-form, or a
knot invariant.  To check this, it suffices by Stokes' Theorem to
examine
the pushforward $\pi_{*}$ along the codimension one faces of
$\FM{2n}{\R}$.  If the integrals vanish on every such face,
\eqref{E:IntroPushforward} yields a knot invariant.
\vskip 5pt
\noindent
This, however, turns out to be too much to hope for.  The
 boundary integrals along some of the faces are nonzero, and one is
next
lead to consider other terms to counter their contribution.  The
correct setting for doing so is provided by first extending the
 chord diagrams to \emph{trivalent diagrams}.  In addition to
 chords on the circle, such diagrams have some number of edges ending
in
 triples at
 vertices off the circle.  The diagrams are oriented by one of the two
 cyclic orientations of the edges emanating from the trivalent vertices
 (see Definition 1.8 in \cite{BN}).

\begin{figure}[h]
\begin{center}
\input{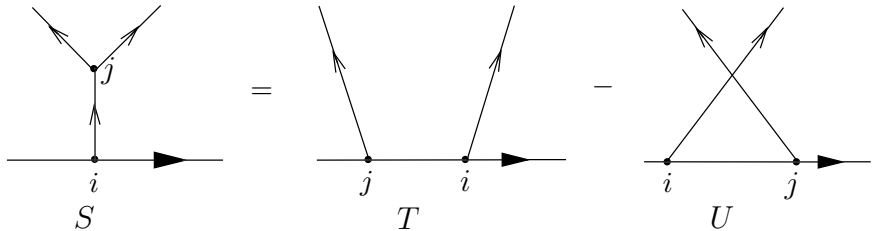}
\caption{STU relation}\label{F:STU}
\end{center}
\end{figure}

 Let $TD_{n}$ be the set of trivalent diagrams with $2n$ vertices.
 Also consider the $STU$ relation in Figure \ref{F:STU}.  Bar-Natan
  uses the fact that the $4T$ relation looks like
 a difference of two such relations to prove
 \begin{thm}[\cite{BN}, Theorem 6.]
 The real vector space generated by trivalent diagrams with $2n$
 vertices modulo the
 $STU$ relation is isomorphic to $\cd{n}$.
 \end{thm}
 The weight systems now extend uniquely to trivalent diagrams modulo
 $STU$.

Given a trivalent diagram with $2n$ vertices, $k$ of which are on
the circle and $s$ are trivalent, Bott and Taubes now construct a
$(k+3s)$-dimensional space $ F[k,s;\K,\R^{3}] $ which fibers over
$\K$ with fiber a configuration space of $2n$ points in $\R^{3}$
with the condition that $k$ of these points are constrained to lie
on a given knot $K$.  Fulton-MacPherson compactification may be
carried out for these spaces as well.  Since trivalent diagrams
are generalizations of chord diagrams, $F[k,s;\K,\R^{3}]$ reduces
to $\FM{2n}{S^{1}}\times \K$ when $s=0$.

So now trivalent diagrams $D$ become prescriptions for pulling back
2-forms from products of spheres.  Each chord and edge of $D$ again
determines a map to $S^{2}$ via the normalized difference of the
appropriate two points in $F[k,s;\K,\R^{3}]$, and each map is used
to pull back a volume form from a sphere.   As described in \cite{BT, Vo2}, the resulting form pulled back by the
product of
all the maps extends smoothly to the
codimension one strata of $F[k,s;\K,\R^{3}]$.  Since $D$ has
$(k+3s)/2$ chords and edges, the resulting form on
$F[k,s;\K,\R^{3}]$
is $(k+3s)$-dimensional.  Its pushforward  thus again yields a
function on $\K$.

Let $I(D,K)$ denote the pullback followed by pushforward as before,
only now $D$ could also be a trivalent diagram.  The following
theorem was proved by Altschuler and Freidel \cite{Alt}, while a different
partial proof was given by D. Thurston \cite{Th}
following the work of Bar-Natan and Bott-Taubes (see also \cite{Vo2}).

\begin{thm}\label{T:UniversalInvariant}  Let $W$ be a weight system in $\W_{n}$.  The sum of
$I(D,K)$
with coefficients $W(D)$, taken over all trivalent diagrams $D$ with
$2n$
vertices
and with a
certain correction term for each $D$, is a knot
invariant.
Further, when $I(D,K)$ is extended to sums of resolutions
of singular knots, it gives a universal type $n$ invariant.
\end{thm}
The proof is a
combination of arguments showing that the integrals along the
various faces of $F[k,s;\K,\R^{3}]$ either vanish, cancel due to the
relations imposed on trivalent diagrams, or can be compensated for by
other integrals.  Next we will see how this theorem extends to the
stages of the Taylor tower for $\K$.

    \section{Taylor tower for the space of knots}\label{S:TaylorTower}

Building on general theory of calculus of functors developed by
Goodwillie \cite{CalcI, CalcII, CalcIII}, Weiss develops in \cite{We} (also see \cite{GKW})
a certain tower
for studying the space of embeddings of one manifold in another.  The
stages of the tower in some sense approximate the original space of
embeddings as long as the codimension of one manifold in the other is
at least three.  This is of course not true in our situation, but the
tower for $\K$ can still be constructed as outlined below.
More on this construction for the particular case of spaces of knots
can be found in \cite{Dev, Vo2}.

\vskip 5pt
\noindent
Let $\{A_{i}\}$, $1\!\leq\! i\!\leq r$, be
a collection of disjoint
closed subintervals of $I$ not containing the endpoints, and define
spaces of
``punctured knots''
$$
E_{S}= \mbox{Emb}(I\!-\! \bigcup_{i\in S}A_{i},S^{3})
$$
for each nonempty subset $S$ of $\{1, \ldots, r\}$.  Here we as
usual mean that the endpoints of the interval are sent to
some basepoint in $S^{3}$ with the same derivative while the other
pieces of $I$ are embedded.

These spaces can be arranged in a subcubical diagram (like a
cubical diagram, but missing a space which maps to all others),
which we will call $EC_r$, since every punctured knot restricts to
a knot with more punctures.  Thus, for example, if $r=3$, we have
a 3-subcubical diagram $EC_{3}$
$$
\xymatrix@=20pt{
          &           &   E_{\{1\}} \ar'[d][dd]
          \ar[dr]  &                  \\
          &  E_{\{2\}} \ar[rr] \ar[dd]  &
          & E_{\{1,2\}}
          \ar[dd] \\
E_{\{3\}} \ar'[r][rr] \ar[dr] &        &   E_{\{1,3\}}
\ar[dr] &                   \\
         &   E_{\{2,3\}} \ar[rr]      &                    &
E_{\{1,2,3\}}
}
$$
Each square side of such a diagram commutes since the order in which
the subintervals of $I$ are removed does not matter.

Now let $\Delta^{S}$ be the face of $\Delta^{r-1}$ for which those
barycentric coordinates indexed by elements of $\{1,\ldots, r\}$
which are not in the subset $S$ are 0.  Recall that the \emph{homotopy
limit} of $EC_r$, which we will denote by $\HO_{r-1}$, can then be defined as the subspace of
$$
\prod_{\emptyset\neq S\subseteq\{1, \ldots, r\}} Maps(\Delta^{S},
E_{S})
$$
consisting of collections of maps $\{\alpha_{S}\}$ such that, for
every map $E_{S}\to E_{S\cup \{i\}}$ in the diagram, the square
$$\xymatrix{
\Delta^{S}\ar[r]^{\alpha_{S}} \ar@{^{(}->}[d] & E_{S} \ar[d] \\
\Delta^{S\cup\{i\}} \ar[r]^{\alpha_{S\cup \{i\}}}   &   E_{S\cup \{i\}}
}
$$
commutes. (More on homotopy limits of diagrams of spaces can be found
in \cite{BK, CalcII}.)

Two observations are immediate:  Any knot can be restricted to a
knot with punctures, so that $\K$ maps to the subcubical diagram
(in fact, $\K$ is also the actual limit of the subcubical diagram,
as long as $r>2$). Consequently, there are canonical maps $
\K\longrightarrow\HO_{r}. $ Also, the subcubical diagram for knots
with up to $r$ punctures is contained in the diagram for knots
with up to $r+1$ punctures, so that there are also projections $
\HO_{r}\longrightarrow \HO_{r-1}. $
\begin{definition} The \emph{Taylor tower for $\K$} is the collection of maps
\begin{equation}\label{E:GoodwillieTower}
\xymatrix{ \HO_{1} & \cdots\ar[l] &  \HO_{r-1}\ar[l] &
\HO_{r}\ar[l] & \cdots\ar[l] }.
\end{equation}
Spaces $\HO_{r}$ will be called the \emph{stages} of the tower.
\end{definition}

\begin{rem}
In the most general setting, the Taylor tower is constructed from
embeddings of codimension zero submanifolds of one manifold in
another.  As mentioned earlier, if the codimension is at least three,
Goodwillie, Klein, and Weiss \cite{GK, GW} show that the maps from the
original space of embeddings to stages of the tower induce more
isomorphisms on homotopy and homology groups the higher one goes in
the tower.  Consequently, the inverse limit of the tower is weakly
equivalent to the original space of embeddings.
In particular, this fact has been used for studying the homology and
homotopy of spaces of knots in $\R^{m}$, $m>3$ \cite{Pasc, SS, Dev, Vo2}.
\end{rem}

In \cite{Vo2}, we show how Bott-Taubes integrals can be extended
to the stages.  We integrate over certain
spaces which generalize $F[k,s;\K,\R^{3}]$. The
$k$ points no longer represent a configuration on a single knot, but
rather on a family of punctured knots $h$, i.e. an element of $\HO_{2n}$.  Denoting by
$I(D,h)$ the
analog of \eqref{E:IntroPushforward}, we prove

\begin{thm}\cite[Theorem 4.1]{Vo2}\label{T:B-TMainTheorem}  Let $W$ be a weight system in
$\W_{n}$ and let $D$ be a
trivalent diagram with $2n$ vertices.
Then the sum of $I(D,h)$ over all $D$, with coefficients $W(D)$, is
an invariant of $\HO_{2n}$, provided a correction term is given for
each $D$.  In particular, this restricts to the invariant of $\K$
described in \refT{UniversalInvariant}
when a point in $\HO_{2n}$ comes from a knot.
\end{thm}

%

The next goal is to see how our invariants $T(W)$ behave when
singular knots, or rather sums of their resolutions, are mapped into the
stages $\HO_{2n}$.  To that end, we introduce a cosimplicial model for the Taylor tower
because it is suitable for our computational purposes.

\section{A Cosimplicial Model for the Taylor
Tower}\label{S:CosimplicialModel}

\subsection{A cosimplicial space of configurations in a manifold}\label{S:CosimplSpace}

Let $M$ be a simply-connected manifold of dimension at least 3,
embedded in $\R^{d}$ for some large
$d$.  We will later specialize to $M=S^{3}$.  Let
\begin{equation}
\begin{array}{l}
p=(p_{0}, p_{1}, \ldots, p_{r}, p_{r+1}), \ \ \  p_{i}\in M
\notag \\
v=(v_{ij}: \ 0\leq i\leq j\leq r+1), \ \ \  v_{ij}\in V, \
|v_{ij}|=1. \notag
\end{array}
\end{equation}

Let $b$ be some fixed point
 in $M$ and $\beta$  a fixed unit vector in $T_{b}M$, the tangent space of $M$
 at $b$.  Define
$$Y^{r}\subset M^{r+2}\times (S^{d-1})^{r+2\choose 2}$$
with subspace topology as follows:  $(p,v)$ is in $Y^{r}$ if it satisfies
\begin{enumerate}
\item $p_{0}=p_{r+1}=b$.
\item $v_{ii}\in T_{p_{i}}M$, $v_{00}=v_{(r+1)(r+1)}=\beta$.
\item $p_{i}=p_{j} \Longrightarrow
\begin{cases}
p_{i}=p_{k}, \ \  v_{ii}=v_{kl}=v_{jj}, \ \   i\leq k\leq l\leq j;
\\[4pt] \text{or}   \\[4pt]
p_{0}=p_{1}=\cdots=p_{i},\ \  v_{00}=v_{kl}=v_{ii}, \ \ 0\leq k\leq l\leq i
\\[-5pt] \text{and}   \\[-5pt]
p_{j}=p_{j+1}=\cdots=p_{r+1},\ \  v_{jj}=v_{mn}=v_{(r+1)(r+1)}, \ \ j\leq m\leq
n\leq r+1.
\end{cases}
$
\end{enumerate}
Let $U^{r}=\{(p,v)\}\in M^{r}\times (S^{d-1})^{r}$ with
\begin{equation}
\begin{array}{l}
p=(p_{1}, p_{2}, \ldots, p_{r}), \ \ \ p_{i}\neq
p_{j}, \  p_{i}\neq b,  \ \mbox{for all} \ i, j \notag \\
v=(v_{1}, v_{2}, \ldots, v_{r}), \ \ \ \ v_{i}\in T_{p_{i}}M, \
|v_{i}|=1 \notag.
\end{array}
\end{equation}
$U^{r}$ is thus the space of configurations of $r$ distinct
points in $M\setminus\{b\}$ labeled with tangent vectors.  In
particular, if the dimension of $M$ is $d$,
$$
U^{r}=F(r, M\setminus\{b\})\times (S^{d-1})^{r}
$$
This space maps to $Y^{r}$ by
$$(p_{1}, \ldots, p_{r}, v_{1}, \ldots, v_{r})\longmapsto (b, x_{1},
\ldots, x_{r}, b, v_{ij})$$
where
\begin{equation}
v_{ij}=
\begin{cases}
v_{i},      &\text{if $0<i=j<m+1$;}\\
\beta,         &\text{if $i=j=0$ or $i=j=m+1$;}\\
\frac{p_{j}-p_{i}}{|p_{j}-p_{i}|}, &\text{if $i<j$.}
\end{cases}
\notag
\end{equation}
It is easy to see that this map is one-to-one and a homeomorphism onto
its image, so
$U^{r}$ can be identified with a subspace of $Y^{r}$.

\begin{definition}\label{D:CosimplicialEntries} Let
$X^{0}=(b, b, \beta)$ and define $X^{r}$ for all
$r>0$ to be the closure of $U^{r}$ in $Y^{r}$.
\end{definition}

We now define the coface and codegeneracy maps,
$\partial^{i}$ and $s^{i}$, which will be given by doubling and
forgetting points and vectors.

Let
$$\partial^{i}\colon Y^{r}\longrightarrow Y^{r+1}$$
be given by
\begin{equation}\label{E:Codegeneracy}
(p_{1},\ldots,p_{i},\ldots,p_{r}, v_{ij})\longmapsto
(p_{0},\ldots, p_{i},p_{i},\ldots,p_{r+1}, v_{ij}'),
\end{equation}
where the $v_{ij}$ map to $v_{ij}'$ by
\begin{equation}
\left(
\begin{smallmatrix}
  &            &         & \ddots & \vdots & \vdots      & \vdots & \vdots \\
  &            &         &        & v_{ii} & v_{i(i+1)}  & \ldots & v_{i(r+1)} \\
  &            &         &        &        & \ddots      & \vdots &   \vdots 
\end{smallmatrix}
\right)
\longmapsto
\left(
\begin{smallmatrix}
      &        &         & \ddots & \vdots & \vdots      & \vdots   & \vdots & \vdots \\
      &        &         &        & v_{ii} & v_{ii}      &v_{i(i+1)}& \ldots & v_{i(r+1)} \\
      &        &         &        &        & v_{ii}      &v_{i(i+1)}& \ldots & v_{i(r+1)} \\
      &        &         &        &        &             & \ddots   & \vdots & \vdots 
\end{smallmatrix}
\right)
\end{equation}
Here one should keep in mind that $p_{0}=p_{r+1}=b$,
$v_{00}=v_{(r+1)(r+1)}=\beta$ as before, and also that
$v_{(i+1)(i+1)}=v_{i(i+1)}=v_{ii}$ in the image array.  The following
is immediate from the definitions:

\begin{lemma} $\partial^{i}(X^{r})\subset X^{r+1}$.
\end{lemma}

This lemma is true for $x_{i}=x_{0}=b$ or $x_{i}=x_{r+1}=b$, so that we
indeed get all the doubling maps we want:
$$
\partial^{i}\colon X^{r}\longrightarrow
X^{r+1}, \ \ \ 0\leq i\leq r+1.
$$
For $s^{i}$, the forgetting maps, the situation is much simpler.
Define
$$
s^{i}\colon X^{r+1}\longrightarrow X^{r}, \ \ \ 1\leq i\leq r
$$
by omitting the entry $p_{i}$ as well as all vectors for which at
least one index is $i$, and then relabeling the result.  This is
clearly continuous.

\begin{prop}  Let $\Xdot$ be the sequence of spaces $X^{0}, X^{1}, X^{2}, \ldots$,
together with doubling and forgetting maps $\partial^{i}$ and $s^{i}$
as given above.  Then $\Xdot$ is a cosimplicial space.
\end{prop}

The proof is a straightforward check of the cosimplicial
identities and is left to the reader.

\subsection{Equivalence of totalizations and stages of the Taylor
tower}\label{S:TotIsHolim}

Recall that $\deldot$ is the cosimplicial space made of closed
simplices $\Delta^{r}$, with cofaces and
codegeneracies the inclusions of and projections onto faces.

Also recall that the \emph{totalization $Tot \Xdot$} of a cosimplicial space
$\Xdot$ is
a subspace of the space of maps from $\deldot$ to $\Xdot$, and can be
defined as
$$\lim_{\longleftarrow}(Tot^{1}\Xdot \longleftarrow
Tot^{2}\Xdot\longleftarrow \cdots).$$
Here $Tot^{r}\Xdot$ is the \emph{$r$th partial
totalization}, or the subspace of
$$\prod_{0\leq i\leq r}Maps(\Delta^{i}, X^{i}), \ \ \ 0\leq i\leq r,$$
determined by the compatibility condition that the squares
\begin{equation}\label{E:TotSquares}
\xymatrix{\Delta^{i} \ar[r]^{\partial^{j}} \ar[d] & \Delta^{i+1}\ar[d]
\\
X^{i} \ar[r]^{\partial^{j}} & X^{i+1}
}\ \ \ \ \
\mbox{and} \ \ \ \ \
\xymatrix{\Delta^{i+1} \ar[r]^{s^{k}} \ar[d] & \Delta^{i}\ar[d]
\\
X^{i+1} \ar[r]^{s^{k}} & X^{i}
}
\end{equation}
commute for all cofaces $\partial^{j}$ and codegeneracies $s^{k}$.
\vskip 5pt
\noindent
We now consider the special case $M=S^{3}$.  Recall that
$K\in \K$ is a map of $I$ in $S^{3}$ embedding the interior and
sending 0 and 1 to a
fixed point $b\in S^{3}$ with the same fixed tangent vectors
$K'(0)=K'(1)=\beta$.  Let $\Delta^{r}$ be parametrized by
$
(0\!=\!x_{0}, x_{1}, \ldots, x_{r}, x_{r+1}\!=\!1).
$

\begin{definition}\label{D:EvaluationMap}
Given a knot $K$, define
$$ev_{r}(K)\colon \Delta^{r} \longrightarrow X^{r}$$
by $$ev_{r}(K)(x_{0}, x_{1}, \ldots, x_{r}, x_{r+1})=
(b, K(x_{1}), \ldots, K(x_{r}), b,v_{ij})$$
where
$$v_{ij}=\begin{cases}
\frac{K'(x_{i})}{|K'(x_{i})|},   &\text{if $0<i=j<r+1$;}\\[8pt]
\frac{K(x_{j})-K(x_{i})}{|K(x_{j})-K(x_{i})|},  &
\text{if $0\leq i\neq j<r+1$\ \ or \  $0<i\neq j\leq r+1$.}
\end{cases}$$

\noindent
Let $$ev_{r}\colon \K \longrightarrow Maps(\Delta^{r},X^{r})$$ be
defined by sending $K$ to $ev_{r}(K)$.
\end{definition}
It is clear that $ev_{r}(K)$ is a continuous map to $X^{r}$ since
it is continuous into $Y^{r}$ and it takes $int(\Delta^{r})$ into
$U^{r}$.  The following is a straightforward check of the
compatibility of $ev_{r}$ with the cosimplicial maps:

\begin{prop}\label{P:evMapstoTot}
Denote by $ev_{[r]}$ the collection of maps $ev_{l}$,
$l\leq r$.  Then $ev_{[r]}$ maps $\K$ to $Tot^{r}\Xdot.$
\end{prop}

We may therefore arrange the maps between partial totalizations of
$\Xdot$ and the maps they admit from $\K$ into a tower much like
the Taylor tower in \eqref{E:GoodwillieTower}, except with
$\HO_{r}$ replaced by $Tot^{r}\Xdot.$ In fact, the next statement
says that $\Xdot$ is a good substitute for the Taylor tower.

\begin{thm}\cite[Theorem 6.1]{Dev}\label{T:Tot=holim}
The $r$th stage $\HO_{r}$ of the Taylor tower for $\K$
is weakly equivalent to $Tot^{r}\Xdot$ for all $r>0$.
\end{thm}

The cosimplicial space Sinha uses in \cite{Dev} is slightly
different than ours, with $X^r$ defined as the closure of the
image of the configuration space of $r$ points in $S^3$ under the
map $i\times \pi_{ij}$ (see Definition
\ref{D:DevCompactification}).  Cofaces and codegeneracies are
doubling and forgetting maps like ours.  Using
\refC{TotSpaces=Configurations} and \cite[Theorem 4.2]{Dev}, it is
easy to see that the two cosimplicial spaces are equivalent.
Namely, they both up to consist of configuration spaces labeled
with tangent vectors up to homotopy, and there are obvious
equivalences between those spaces respecting the cosimplicial
maps.  The advantage of our definition is that it is made with
knots and evaluation maps in mind (points moving along a knot can
only collide in one direction).

We next prove in \refP{PuncturedKnots=Configurations} and
\refC{TotSpaces=Configurations} that $E_{S}$ (spaces of punctured
knots) and $X^{r}$ are configurations labeled with tangent
vectors. (The above theorem then in effect says that these
configuration spaces are ``put together the same way'' in the
homotopy limits and partial totalizations. The main observations
and tools in Sinha's proof are that the restriction maps between
punctured knots look like doubling maps up to homotopy and that a
truncated cosimplicial space can be redrawn as a subcubical
diagram whose homotopy limit is equivalent to the original
totalization.  Sinha thus constructs equivalences going through an
auxiliary tower whose stages are homotopy limits of cubical
diagrams of compactified configuration spaces with doubling maps
between them.)

\vskip 5pt
\noindent
Recall that $E_{S}$ is the space of embeddings of the complement
of $s=|S|$ closed
subintervals in
$S^{3}$, where $S$ is a nonempty subset of $\{1, \ldots,r\}$.  In
other words, each point in $E_{S}$ is an embedding of $s+1$ open subintervals
(or half-open, in case of the first and last subinterval), which we
index in order from 1 to $s+1$.

For each $S$ and
each open subinterval, choose points $x_{i}$, $1\leq i\leq s-1$, in
the interiors.  The choices should be compatible, namely if $S\subset T$,
then $\{x_{i}\}_{i\in S}\subset \{x_{i}\}_{i\in T}$.  Also let $x_{0}$
and $x_{s}$ be the left and right endpoints of $I$, respectively.

\begin{prop}\label{P:PuncturedKnots=Configurations}
$E_{S}\simeq F(s-1,S^{3}\!\setminus\!\{b\})\times (S^{2})^{s-1}.$
\end{prop}

\begin{proof}
To specify an embedding at a point, it is sufficient and necessary to
specify the image of the point as well as the nonzero derivative of the
embedding at that point.  For more points, the images additionally
must be distinct.  Letting $e$ be an embedding in $E_{S}$, we thus
have an evaluation map
\begin{align*}
ev_{S}\colon  E_{S} & \longrightarrow
F(s-1,S^{3}\!\setminus\!\{b\})\times (S^{2})^{s-1} \\
e & \longmapsto
\left(e(x_{0}),\ldots, e(x_{s}), \frac{e'(x_{0})}{|e'(x_{0})|}
\ldots, \frac{e'(x_{s})}{|e'(x_{s})|}\right).
\end{align*}
Here we need to remember that the endpoints always map to the
origin with the same tangent vector.  (As a matter of fact, we could have omitted
either one of the endpoints from the above map.)

To see that $ev_{S}$ is an equivalence, observe that its fiber consists of
all embeddings in $E_{S}$ going through $s-1$
distinct points in $S^{3}$ with some specified tangent vectors.  But
this is a contractible space as any such embedding can be ``shrunk
back'' to the specified points while preserving the tangent vectors at
those points.
\end{proof}

We next prove an almost identical statement for the spaces in $\Xdot$ as a
corollary of
\begin{prop}\label{P:weak}
The inclusion ${\displaystyle U^{s}\hookrightarrow X^{s}}$ induces a
weak homotopy equivalence.
\end{prop}
Since $U^{s}$ is a space of $s$ configuration points in the interior
of $S^{3}\!\setminus\!\{b\}$ with vectors attached, we deduce the following:
\begin{cor}\label{C:TotSpaces=Configurations} ${\displaystyle X^{s}\sim
F(s, S^{3}\!\setminus\!\{b\})\times (S^{2})^{s}.}$
\end{cor}
To prove Proposition~\ref{P:weak}, we will need a
technical result, whose proof is straightforward:
\begin{lemma}\label{L:weak} Let $U$ be an open subset of a space $X$.  Suppose
$$H\colon X\times [0,1]\longrightarrow X$$
is a map such that, for each $p\in X$, $H(p,0)=p$ and there is a neighborhood
$V_{p}$ of $p$ and $\epsilon_{p}\in (0,1]$ with
$$H(y,t)\in U \ \ \ \text{for all} \ \ \ (y,t)\in
(V_{p},(0,\epsilon_{p}]).$$
Then the inclusion $U\hookrightarrow X$ induces a weak equivalence.
\end{lemma}

\begin{proof}[Proof of Proposition~\ref{P:weak}.]
We wish to define a homotopy $$H\colon X^{s}\times [0,1]\longmapsto
X^{s}$$ satisfying the conditions of the previous Lemma, namely that
each point $p\in X^{s}$ has a neighborhood $V_{p}$ on which $H$ is a homotopy
into $U^{s}$ for some subinterval $[0,\epsilon_{p}]$.  To do this, we
use the vectors $v_{ii}$ to ``separate'' points $p_{i}$ and $p_{i+1}$
when $p_{i}=p_{i+1}$,
thereby mapping them from the boundary of $X^{s}$ to its interior,
$U^{s}$.

We first do this locally.  For each $p=(p_{0}, \ldots, p_{s+1})$ let
$$w_{i}=\sum_{i=0}^{s+1}c_{i}(p)v_{ii},$$
where each $c_{i}(p)$ assigns a real number to point $p_{i}$, and they
together
 satisfy the conditions
\begin{equation}\label{E:HomotopyConditions}
(c_{i+1}-c_{i})(p)>0 \text{ \ if\  } p_{i}=p_{i+1}, \ \ \text{ and } \ \
c_{0}=c_{s+1}=0.
\end{equation}
Notice that, if $p_{j}\neq p_{j+1}$, $(p,v)$ always has a neighborhood
in $X^{s}$
such that for all $(y,v)$ in that neighborhood, $y_{j}\neq
y_{j+1}$.  We therefore also have a neighborhood $V_{p}$ of $(p,v)$ on
which $H$, defined as
$$H(p_{i}, v_{ij}, t)=(p_{i}+tw_{i}, v_{ij}),$$
is a homotopy into $U^{s}$ for at least some time interval
$[0,\epsilon_{p}]$.  Further, $c_{i}$ can then be chosen to be constant
on $V_{p}$, and
if $p_{j}\neq p_{j+1}$, $c_{j}$ may be set to be 0 so that $H$ is
continuous on $V_{p}$.

We may now cover $X^{s}$ with open sets $\{V_{\alpha}\}$ determined by
$V_{p}$, and use a partition of unity to define a continuous
$H$ on all of $X^{s}$.  It is immediate that the conditions
\eqref{E:HomotopyConditions} are satisfied on the intersections of
the $V_{\alpha}$.
\end{proof}

\subsection{Tangential data and framed knots}\label{S:Framing}

To make future arguments work out easier, in this section
we
remove the tangential data appearing in \refP{PuncturedKnots=Configurations} and
\refC{TotSpaces=Configurations}.  The price we will have to pay is
that the space of ordinary knots will
become a space which on $\pi_{0}$ looks like the space of framed
knots whose framing number is even.
This is why we worked with framed
knots throughout the previous sections, and, as
explained earlier, the
class of chord diagrams was larger for us than it would have been
had the framing not been considered.
\vskip 5pt
\noindent
First note that \S\ref{S:TaylorTower} could be repeated
with ``immersions'' instead of ``embeddings.''  So let

\begin{itemize}
\item $Im$ be the
space of immersions of $I$ in $S^{3}$ with usual conditions on
endpoints,

\item $Im_{S}$ the
space of immersions of $I$ with $s=|S|$ subintervals removed,

\item $ImC_{r}$
the subcubical diagram obtained by considering $Im_{S}$ for all nonempty
subsets $S$ of
$\{1,\ldots,r+1\}$ with restriction maps between them, and

\item $holim(ImC_{r})$ its
homotopy limit.

\end{itemize}

Recall that $Im$ is homotopy equivalent to $\Omega S^{2}$ since
immersions are determined by their derivatives.  Similarly, in analogy with Proposition
\ref{P:PuncturedKnots=Configurations}, we have
\begin{equation}\label{E:Immersions=Spheres}
Im_{S}\simeq (S^{2})^{s-1}.
\end{equation}
The configuration space is no longer present in the equivalence because
of the lack of the
injectivity condition for immersions.

Before we prove two useful statements, we need
\begin{prop}[\cite{CalcII}, Proposition 1.6.]\label{P:StronglyCartesian}
Suppose
$X_{\emptyset}$ completes a subcubical diagram $C_{r}$ of spaces $X_{S}$
indexed by nonempty subsets of $\{1, \ldots, r\}$ into a cubical diagram.
If, for every $S$,
$$
X_{S}\longrightarrow holim(X_{S\cup\{j\}}\to X_{S\cup\{i,j\}}\leftarrow
X_{S\cup\{i\}})
$$
is a weak equivalence, then
$
X_{\emptyset}\to holim(C_{r})
$
is a weak equivalence as well.
\end{prop}

\begin{prop}\label{P:ImmersionsHolim}
$Im\longrightarrow holim(ImC_{r})$ is a weak equivalence.
\end{prop}

\begin{proof}
Unlike for spaces of punctured embeddings, $Im_{S}$ is the limit of
\begin{equation}\label{E:ImmersionSquare}
\xymatrix{
 Im_{S\cup\{j\}}\ar[r]      &  Im_{S\cup\{i,j\}} &  Im_{S\cup\{i\}} \ar[l]
}
\end{equation}
Further, the map $Im_{S}\to
Im_{S\cup\{i\}}$ is a fibration for all $i$.  It follows that $Im_{S}$ is
weakly equivalent to the homotopy limit of diagram
\eqref{E:ImmersionSquare}.
This holds for any square in $ImC_{r}$.  In addition, we can
complete the subcube $ImC_{r}$ by adding $Im$ as the initial space.
For the extra square diagrams
produced this way, we again have that $Im$ is equivalent to their
homotopy limits.
\refP{StronglyCartesian} then finishes the proof.
\end{proof}

\begin{prop}  There is a map of homotopy limits
$$
holim(EC_{r})\longrightarrow holim(ImC_{r}),
$$
whose homotopy fiber is the homotopy limit of a subcubical diagram
$FC_{r}$ of spaces which are homotopy equivalent to
configuration spaces of up to $r$ labeled points in $S^{3}\!\setminus\!\{b\}
$.
\end{prop}

\begin{proof}
Since every embedding is an immersion, there are inclusions
$E_{S}\longrightarrow Im_{S}$
for each $S$.  By Proposition
\ref{P:PuncturedKnots=Configurations} and equation
\eqref{E:Immersions=Spheres} this map is,
up to homotopy, the projection
$$
F(s,S^{3}\!\setminus\!\{b\})\times (S^{2})^{s}\longrightarrow (S^{2})^{s},
$$
with fiber the configuration space $F(s,S^{3}\!\setminus\!\{b\})$.
However, for each $S$ such that $|S|=s$, we get a different
configuration space of \emph{labeled} points in $S^{3}\!\setminus\!\{b\}$.  It follows
that there is a map of subcubes,
\begin{equation}\label{E:CubesMap}
EC_{r}\longrightarrow ImC_{r},
\end{equation}
with homotopy fiber a subcubical diagram $FC_{r}$ of configuration spaces
as desired.
\end{proof}

If we denote by $E$ the space of unframed knots and by $\overline{\K}$ the homotopy fiber of the projection
$E\to Im$, the previous two propositions can be summarized in the diagram
\begin{equation}
\xymatrix@R=30pt{
\overline{\K}  \ar[r] \ar[d]  &  E  \ar[r]  \ar[d]  & Im\simeq\Omega S^{2} \ar[d]^{\sim}  \\
holim(FC_{r})  \ar[r]         &  holim(EC_{r}) \ar[r] & holim(ImC_{r})
}
\end{equation}
where the left spaces are the homotopy fibers of the two right
horizontal maps.

The middle column in the diagram represents the spaces and maps we
would have originally considered, but since we wish to remove the tangent
spheres,
we now
shift our attention to the left column.  Since the map $E\to Im$ is
null-homotopic \cite[Proposition 5.1]{Dev3}, we have
$$
\overline\K \simeq E\times \Omega Im.
$$
\begin{prop}\label{P:Framing}
$\pi_{0}(\overline{\K})$ may be identified with the isotopy classes
of framed knots whose framing number is even.
\end{prop}

\begin{proof} The space of framed knots up to isotopy is
$\pi_{0}(E)\times \Z$
since every knot can be framed to yield any framing number and this
number is an isotopy invariant.
To prove the proposition, it hence suffices to exhibit a
bijection
\begin{equation}
\xymatrix{
\pi_{0}(\overline{\K}) \ar[r]^(0.4){f} & \pi_{0}(E)\times 2\Z.
}
\end{equation}
A point in $\overline{\K}$ consists of an embedding
$\widetilde{K}$ along
with a path $\alpha$ from $\widetilde{K}$
through immersions to the
basepoint in $Im$, namely the unknot.  This path produces a possibly
singular surface where the singularities come from $\widetilde{K}$ passing
through itself.  The standard framing of the
unknot (a copy of the unknot, displaced slightly in one direction so that the
framing number of the unknot is 0) can be used to obtain a
framing of $\widetilde{K}$ by dragging the
framing of the unknot along $\alpha$.  However, as $\alpha$ goes through a
singularity, the framing number of a knot changes by 2.  It follows that
$\widetilde{K}$ must have an even framing number.

Next we
have an exact sequence of homotopy groups whose end is
$$\cdots\to\Z\to\pi_{0}(\overline{\K})\to
\pi_{0}(E)\to 0.$$
Here we have used
$
\pi_{0}(\Omega S^{2})=0$ and $\pi_{1}(\Omega S^{2})=\Z.$

Now $\pi_{1}(Im)=\Z$ acts on
$\pi_{0}(\overline{\K})$ and it does so in such a way that the action
of $1$ adds 2 to the framing number, while the action of $-1$ subtracts 2. This is
because the generator for $\pi_{1}(Im)$ is the loop of the unknot that
introduces a twist (the first Reidemeister move), passes the crossing
obtained that way through
itself, and then untwists the result having changed the framing
number by 2.  Depending on the orientation of
the unknot and on the direction of the twist (right or
left-handed), this loop corresponds to one of the generators of $\Z$.

Given an element $(K,2n)$ in $\pi_{0}(E)\times 2\Z$,
where $K$ represents a knot type, there
is a class represented by
$(K,\alpha)$ in $\pi_{0}(\overline{\K})$ with $\alpha$ a path of $K$ to the unknot
through immersions.  The path then provides some framing of $K$ with
the framing number $2m$ so that the action of $n-m$ on $(K,\alpha)$,
composed with $f$, gives $(K,2n)$.  So $f$ is surjective.

Now suppose two representatives of classes in $\pi_{0}(\overline{\K})$,
$(\widetilde{K}_{1},\alpha_{1})$ and
$(\widetilde{K}_{2},\alpha_{2})$, give isotopic knots $K_{1}$ and
$K_{2}$ with the same even framing number.  By exactness of the
homotopy sequence, there exists an integer which acts on
$(\widetilde{K}_{1},\alpha_{1})$ to produce
$(\widetilde{K}_{2},\alpha_{2})$.  But since $K_{1}$ and
$K_{2}$ have the same framing number, this integer must be 0, so that
$(\widetilde{K}_{1},\alpha_{1})$ and
$(\widetilde{K}_{2},\alpha_{2})$ represent the same element of
$\pi_{0}(\widetilde{E})$.  Thus $f$
is injective.
\end{proof}
To set the notation for the next section, let
$$
T_{r}=holim(FC_{r+1}),
$$
We thus have the modified Taylor tower for the space of knots,

\begin{equation}\label{E:ModifiedTower}
\xymatrix{
  T_{2} &  T_{3} \ar[l] & \cdots \ar[l] &   T_{r} \ar[l] & \cdots  \ar[l]
}
\end{equation}
where each space is equipped with a map from $\overline{\K}$ (a collection
of compatible evaluation maps) making all the resulting triangles commute.

We can now do the same for the cosimplicial space $\Xdot$ and see
that \refT{Tot=holim} is preserved.  Namely, since each $X^{s}$ in
$\Xdot$ is equivalent to $F(s, S^{3}\!\setminus\!\{b\})\times
(S^{2})^{s}$, we consider the cosimplicial space $Y^{\bullet}$
with only $(S^{2})^{s+1}$ as its $s$th space.  The maps are
inclusions and projections. Then $Tot^{r}Y^{\bullet}$ is $\Omega
S^{2}$ for all $r>1$ (since $Y^{\bullet}$ is just the space of
maps from the standard simplicial model for $S^1$ to $S^2$). This
in fact allows us to say that there is a tower
$$
Tot^{2}F(s, S^{3}\!\setminus\!\{b\})^{\bullet}\longleftarrow
Tot^{3}F(s, S^{3}\!\setminus\!\{b\})^{\bullet}\longleftarrow
\cdots \longleftarrow
Tot^{r}F(s, S^{3}\!\setminus\!\{b\})^{\bullet}\longleftarrow
\cdots
$$
where $Tot^{r}F(s, S^{3}\!\setminus\!\{b\})^{\bullet}$ is the homotopy fiber of the map
$$
Tot^{r}\Xdot\longrightarrow Tot^{r}Y^{\bullet}=\Omega S^{2}.
$$
These homotopy fibers are exactly the partial totalizations of a
cosimplicial space consisting of configuration spaces and doubling
maps in fixed directions.

We then have a diagram
$$
\xymatrix{
\HO_{r}  \ar@{^{(}->}[d] \ar[r]^(0.4){\simeq}  &  Tot^{r}\Xdot
\ar@{^{(}->}[d]  \\
holim(ImC_{r})  \ar[r]^(0.6){\simeq}  &   Tot^{r}Y^{\bullet}
}
$$
where the top equivalence is \refT{Tot=holim}, and the bottom comes
from the fact that both spaces are $\Omega S^{2}$ (see
\refP{ImmersionsHolim}).
It follows that the homotopy fibers of the two vertical inclusions
 are weakly equivalent.  These fibers are precisely
$T_{r}$ and $Tot^{r}F(s, S^{3}\!\setminus\!\{b\})^{\bullet}$.

Therefore removing the spheres gives equivalent partial
totalizations and
 homotopy limits of subcubical diagrams.
From now on, we will refer to both as $T_{r}$ and which model we mean
will be clear from
the context.  To simplify notation, we will continue to denote by
$\Xdot$ the new cosimplicial space differing from the old one in that
the tangential data has been removed.

\section{Finite type invariants and the Taylor tower}

Here we associate a spectral sequence to the tower
\eqref{E:ModifiedTower} and deduce all the main results.  This spectral
sequence arises from any
cosimplicial space and this is precisely the reason why we
introduced $\Xdot$.  Its $E_{2}$ term will turn out to be the
connection between the Taylor tower and finite type invariants.

\subsection{Weight systems and the cohomology spectral
sequence}\label{S:Cd'sinSS}

This section briefly recalls the procedure for turning a
cosimplicial space into a double complex from which the $E_{1}$
term of the cohomology spectral sequence is constructed. The
spectral sequence then converges to the cohomology of the total
complex of this double complex, which can be thought of as an
algebraic analog of the totalization $Tot\Xdot$.  In favorable
cases, which may or may not be case here, the cohomology of the
total complex will be isomorphic to the cohomology of the
totalization of the cosimplicial space.  Sinha \cite{Dev} shows
this to be true for the Taylor tower for spaces of knots in
$\R^{m}$, $m>3$, and we use this fact in \cite{Pasc} and
\cite{Vo2}. For more details about the construction of the
spectral sequence for a cosimplicial space, see \cite{Bous,BK}.
\vskip 5pt \noindent Recall that $\Xdot$ consists of spaces
$X^{p}=F(p, S^{3}\!\setminus\!\{b\})$, $p\geq 0$ (tangent vectors
have been removed).  But $S^{3}\!\setminus\!\{b\}$ is homeomorphic
to $\R^{3}$.  To simplify notation, we thus set
$$
F(p)=F(p, \R^{3})\simeq F(p, S^{3}\!\setminus\!\{b\}).
$$
We briefly recall
the ring structure of the cohomology of $F(p)$.  More details can be found in
\cite{Coh}.

\vskip 5pt
\noindent
Let $\nu$ denote the invariant unit volume form on $S^{2}$ as before.
Let $(x_{1}, \ldots, x_{p})\in F(p)$.  Then we have maps
$$
x_{ij}\colon F(p)\longrightarrow S^{2}
$$
given by
$$
(x_{1}, \ldots, x_{p})\longmapsto \frac{x_{j}-x_{i}}{|x_{j}-x_{i}|}.
$$
Each can be used to pull back $\nu$ to $F(p)$.  Let
$\alpha_{ij}=x_{ij}^{*}\nu$.  Then
 $H^{q}(F(p))$ is 0 if $q$ is odd, and it is
generated by products of $\alpha_{ij}$ with relations
\begin{gather}
\alpha_{ij}^{2}=0,\label{E:CohomCond2} \\
\alpha_{ij}=-\alpha_{ji}, \label{E:CohomCond3} \\
\alpha_{ij}\alpha_{kl}=\alpha_{kl}\alpha_{ij}, \label{E:CohomCond4} \\
\alpha_{ij}\alpha_{jk}+\alpha_{jk}\alpha_{ki}+\alpha_{ki}\alpha_{ij}=0,
\label{E:CohomCond5}
\end{gather}
if $q$ is even.

The Poincar\'{e} polynomial for $F(p)$ is
\begin{equation}\label{E:PoincarePoly}
P_{p}(t)=(1+t^{2})(1+2t^{2})\cdots(1+(p-1)t^{2}).
\end{equation}
\vskip 5pt
\noindent
We can now set up a double cochain complex by applying cochains
to the cosimplicial space:
\begin{equation}
C^{*}\Xdot=(C^{*}X^{1} \stackrel{\partial^{h}}{\longleftarrow}
C^{*}X^{2} \stackrel{\partial^{h}}{\longleftarrow}  C^{*}X^{3}
\stackrel{\partial^{h}}{\longleftarrow} \cdots).
\end{equation}
It is not hard to see that
$$
\partial^{h}=(\partial^{0})^{*}-(\partial^{1})^{*}+(\partial^{2})^{*}-\cdots.
$$
gives a differential in the horizontal direction. Let the vertical
differential be denoted by $\partial^{v}$.

Denoting the total complex
of this double complex by $TotC^{*}\Xdot$, there is a second
quadrant spectral sequence
converging to the cohomology of $TotC^{*}\Xdot$.
This cohomology has two filtrations coming from a choice of starting at $E_{1}$ by
taking the cohomology with respect to $\partial^{h}$ or $\partial^{v}$.
We use the latter and set up the spectral sequence so that
$$
E_{1}^{-p,q}=H^{q}(F(p)).
$$
With \eqref{E:PoincarePoly} in hand, the $E_{1}$ term is easily
computable.
However, the double complex $C^{*}\Xdot$ can in some way be ``simplified''
prior
to computing this page.  Namely,
instead of considering all of $C^{*}\Xdot$, we only look at a certain
 subcomplex, $NC^{*}\Xdot$, called the \emph{Dold-Kan
normalization}.  It is constructed in the following way:
\begin{align}
NC^{*}X^{0}= & C^{*}X^{0} \notag \\
NC^{*}X^{p}= & C^{*}X^{p}/G^{p}, \ \ p>0, \label{E:Normalization}
\end{align}
where $G^{p}$ is the subgroup generated by
\begin{equation}
\sum_{i=0}^{p-1}Im(C^{*}X^{p-1}\stackrel{(s^{i})^{*}}
{\longrightarrow}C^{*}X^{p})).
\end{equation}

This normalization can be applied to any double complex obtained from
 a simplicial or a cosimplicial space.  The important feature is that
 $TotNC^{*}\Xdot$ has the same
cohomology as $TotC^{*}\Xdot$.  The normalization in general can be
thought of as throwing away
the
 degenerate part of a
(co)simplicial complex.

Let $G^{p,n}$ be the subgroup of $H^{2n}(F(p))$ generated
by the images of
$$(s^{i})^{*}\colon H^{2n}(F(p-1))\longrightarrow H^{2n}(F(p)).$$
Then the normalized $E_{1}$ has entries
$$
E_{1}^{-p,2n}=H^{2n}_{norm}(F(p))=H^{2n}(F(p))/G^{p,n}.
$$
Choose a basis for $H^{2n}(F(p))$ using
relations \eqref{E:CohomCond2}--\eqref{E:CohomCond5} by letting
$$
\alpha=\alpha_{i_{1}j_{1}}\alpha_{i_{2}j_{2}}\cdots
\alpha_{i_{k}j_{k}}, \ \ i_{m}, j_{m}\in \{1, \ldots, p\},
$$
be a generator of $H^{2n}(F(p))$.  Then $\alpha$ is a basis element if
\begin{gather}
i_{a}=i_{b}\Longrightarrow  j_{a}\neq j_{b},  \label{E:BasisCond2}\\
i_{a}< j_{a},  \label{E:BasisCond3}\\
j_{1}< j_{2}< \cdots < j_{k},  \label{E:BasisCond4}
\end{gather}
The first observation is

\begin{lemma}\label{L:below}  If $2n<p$, then $H^{2n}_{norm}=0$.
\end{lemma}

\begin{proof}
Since $2n<p$, it follows that some
point $x_{r}$ in $F(p)$ was not used in any $x_{ij}$.  Now apply the
codegeneracy $s^{r}\colon F(p)\longrightarrow F(p-1)$
by forgetting $x_{r}$ and relabeling the result.  Then $x_i\in F(p)$
remains $x_i\in F(p-1)$ if $i<r$, and becomes $x_{i-1}\in F(p-1)$ if
$i>r$.

Now let $\alpha'$ be a class in $H^{2n}(F(p-1))$
constructed from the maps $x_{ij}'$ such that
\begin{equation}\label{E:UseCasesAgain}
x_{ij}'=
\begin{cases}
x_{i(j-1)},       &\text{if $i<r<j$;}\\
x_{(i-1)(j-1)},   &\text{if $r<i<j$;}\\
x_{ij},           &\text{if $i<j<r$.}
\end{cases}
\end{equation}
But by construction, composing the product of the $x_{ij}'$ with
$s^r$ yields precisely $\alpha$.  Thus $\alpha$ is in the image of
$$(s^{r})^{*}\colon H^{2n}(F(p-1))\longrightarrow H^{2n}(F(p)).$$
\end{proof}

The discussion in the preceding proof can be extended to other
cases.
Regardless of the relation between $2n$ and $p$, we can still
conclude that
if a class $\alpha$ is obtained without using all the points in
the configuration, some codegeneracy induces
a map on cohomology with $\alpha$ in its image.

On the other
hand, if $\alpha$ is obtained by using all the points, then for every
codegeneracy $s^{i}$,
a map on $F(p)$ involving $x_{i}$ becomes undefined on $F(p-1)$.
Thus $\alpha$ is not in the image of $(s^{i})^{*}$ for any $i$, nor is
it in the subspace generated by them as the image of each codegeneracy
is generated by a part of the basis.  Hence $\alpha$ survives to
$H^{2n}_{norm}(F(p))$.  For the same reason, all such $\alpha$ remain
independent after normalization.
It follows therefore that a basis element $\alpha$ in
$H^{2n}_{norm}(F(p))$
can be described by simply adding one more requirement to conditions
\eqref{E:BasisCond2}--\eqref{E:BasisCond4}:
\begin{equation}\label{E:BasisCond6}
\text{For every $x_{i}\in F(p)$, $i$ must occur as a subscript in
$\alpha$}.
\end{equation}

It is convenient at this stage to switch the point of view from
cohomology to homology for establishing a clear connection between
the Taylor tower (and its cosimplicial
model) and finite type knot invariants.

Given a basis
element $\alpha=\alpha_{i_{1}j_{1}}\alpha_{i_{2}j_{2}}\cdots
\alpha_{i_{n}j_{n}}\in H^{2n}_{norm}(F(p))$, let
$a=a_{i_{1}j_{1}i_{2}j_{2}\cdots i_{n}j_{n}}$
denote its dual in
$H_{2n}^{norm}(F(p))$. (By $H_{2n}^{norm}(F(p))$ we mean the homology
obtained by considering
chains on $\Xdot$, using the alternating sums of $(\partial^{i})_{*}$
to obtain a chain complex in one direction, and then normalizing it via a dual version
of \eqref{E:Normalization}.)
An element $a\in
H_{2n}^{norm}(F(p))$
can be thought of as a chord diagram oriented by $p$
 labeled points which are
connected by oriented chords as prescribed by pairs of indices in $a$.
Thus $i_{1}$ is
connected to
$j_{1}$ and the arrow on the chord points from $i_{1}$;  $i_{2}$ is
connected to $j_{2}$ with the arrow pointing from $i_{2}$, etc.
We also impose the following sign convention:  If two basis elements
$a$ and $a'$ differ by $\tau$ transpositions of their subscripts,
then $a=(-1)^{\tau}a'$.  Thus we choose any element $a$ and a sign for
it, and then apply the sign convention.  It is easy to see that this is a well-defined
assignment.

Conversely, take a free module generated by
chord diagrams
with $p$ labeled vertices and $n$
oriented chords, and consider its quotient obtained
by imposing the relations
\begin{enumerate}
\item If a chord in a diagram $D$ connects the same vertex, then
$D=0$;
\item If two chords in $D$ connect the same vertices, then $D=0$;
\item If two diagrams $D_{1}$ and $D_{2}$ differ by a change of
orientation
of a chord or by a transposition of two vertex labels, then $D=-D'$;
\item Let $i<j<k$.  If
\begin{gather*}\text{$D_{1}$ contains two chords connecting
vertices $i$, $j$ and $j$, $k$;}\\
\text{$D_{2}$ contains two chords connecting
vertices $i$, $j$ and $i$, $k$;}\\
\text{$D_{3}$ contains two chords connecting
vertices $i$, $k$ and $j$, $k$;}
\end{gather*}
then $D_{3}=-D_{1}-D_{2}$.
\end{enumerate}

It follows by construction that $H_{2n}^{norm}(F(p))$ is
isomorphic to this quotient.
In the most relevant case for us, $p=2n$, every vertex of a chord
diagram must then have exactly one chord connecting it to some other
vertex.  But this is precisely the description of $CD_{n}$.  We thus
have the
following simple but
important statement:
\begin{prop}\label{P:CD} $E^{1}_{-2n,2n}=H_{2n}^{norm}(F(2n))$ is
generated by $CD_{n}$.
\end{prop}

The next proposition states that the first differential introduces
precisely the $4T$ relation of \S\ref{S:FinTypeInv's}.  We prove the
dual version because the combinatorics of
the proof will be more transparent.
Remember that
the space of weight systems $\W_{n}$ was defined
as functions on $\R[CD_{n}]$ vanishing on the $4T$ relation.

\begin{prop}\label{P:E2} $E_{2}^{-2n,2n}\cong \W_{n}.$
\end{prop}

\begin{proof}
Since $d^{1}=d^{h}$ and
$H^{2n}_{norm}(F(2n+1))=0$, the image of $d^{1}$
in
 $H^{2n}_{norm}(F(2n))$ is 0 and so
$$E_{2}^{-2n,
2n}=ker(H^{2n}_{norm}(F(2n))\stackrel{d^{1}}{\longrightarrow}H^{2n}_{norm}
(F(2n-1))).$$

Let $\langle \ ,\ \rangle$ denote evaluation of cohomology on
homology, $f$ an element of $H^{2n}_{norm}(F(2n))$, and $B$ the basis
of $H_{2n}^{norm}(F(2n-1)).$  Then
$$
f\in ker(d^{1}) \Longleftrightarrow \langle f, d_{1}a\rangle =0
\ \ \forall a\in B,
$$
where $d_{1}$ is the alternating sum of maps induced by cofaces on
homology.

Recall that $\partial^{r}$
sends $x_{i}\in F(2n-1)$ to $x_{i}\in
F(2n)$ if $i\leq r$ and to $x_{i+1}$ if $i>r$.  It follows that
composing any map
$x_{ij}\colon F(2n)\longrightarrow S^{2}$
with $\partial^r$ gives a map
$$x_{ij}'\colon F(2n-1)\longrightarrow S^{2},$$
where
$
x_{ij}'
$
are given exactly as in \eqref{E:UseCasesAgain}, except $i=r$ and
$j=r$ are also possibilities in the first and last case respectively.

Then
$(\partial^{r})^{*}$
sends a basis element $\alpha$ to a basis element $\alpha'$ by
reindexing each 2-dimensional class $\alpha_{ij}$ in the same
way as in \eqref{E:UseCasesAgain}.

Notice that $\alpha'$ will have exactly one index that is repeated, as
it should.  Namely, indices $r$ and $r+1$ in $\alpha$ will both
become $r$ in $\alpha'$.  However, any time $\alpha'$
has as factors two classes $\alpha_{ir}'$ and $\alpha_{jr}'$,
$i<j<r$, we may
rewrite it in terms of two other basis elements using
\eqref{E:CohomCond5}.


Now let
$$
\alpha'=\alpha_{i_{1}j_{1}}'\alpha_{i_{2}j_{2}}'\cdots
\alpha_{i_{n}j_{n}}'
$$
be a basis element of $H^{2n}_{norm}(F(2n-1))$.  Two of the indices
must be the same---either $\alpha'$ contains
$\alpha_{bc}'\alpha_{bd}'$
or $\alpha_{bc}'\alpha_{cd}'$ (the latter two 2-forms do not
necessarily have to be next to each other in $\alpha'$).  Consider
the first case.  There are four basis elements $\alpha\in
H^{2n}(F(2n))$ which will yield this $\alpha'$ under some
$(\partial^{i})^{*}$:
\begin{align}
\alpha_{1}=& \alpha_{i_{1}j_{1}}\alpha_{i_{2}j_{2}}\cdots
\alpha_{b(c+1)}\alpha_{(b+1)(d+1)}
\cdots \alpha_{i_{n+1}j_{n+1}}, \notag \\
\alpha_{2}=& \alpha_{i_{1}j_{1}}\alpha_{i_{2}j_{2}}\cdots
\alpha_{b(c+1)}\alpha_{(b+1)(d+1)}
\cdots \alpha_{i_{n+1}j_{n+1}}, \label{E:FourElements}\\
\alpha_{3}=& \alpha_{i_{1}j_{1}}\alpha_{i_{2}j_{2}}\cdots
\alpha_{bd}\alpha_{c(d+1)}
\cdots \alpha_{i_{n+1}j_{n+1}}, \notag \\
\alpha_{4}=& \alpha_{i_{1}j_{1}}\alpha_{i_{2}j_{2}}\cdots
\alpha_{bd}\alpha_{c(d+1)}
\cdots \alpha_{i_{n+1}j_{n+1}}. \notag
\end{align}
Here, $\alpha_{1}$ and $\alpha_{2}$ give exactly $\alpha'$ under
$(\partial^{b})^{*}$.  However, $(\partial^{d})^{*}$, applied to
$\alpha_{3}$ and $\alpha_{4}$, produces an element which is identical
to $\alpha'$ except for the factor $\alpha_{bc}'\alpha_{bd}'$.  This
factor is in fact
replaced by $\alpha_{bd}'\alpha_{cd}'$.  But this is not an element
in the basis of $H^{2n}_{norm}(F(2n-1))$, and it can be rewritten as
 $-\alpha'-\alpha''$, where $\alpha''$ contains
$\alpha_{bc}'\alpha_{cd}'$.  Note that the four forms $\alpha_{i}$
only differ in the two factors as indicated.

Similarly, in the second case of a possible repeated index in
$\alpha'$, we again get four basis elements $\alpha_{i}$ with
$\alpha_{3}$ and $\alpha_{4}$ the same as above, and with
\begin{align}
\alpha_{1}= & \alpha_{i_{1}j_{1}}\alpha_{i_{2}j_{2}}\cdots
\alpha_{b(c+1)}\cdots\alpha_{c(d+1)}
\cdots \alpha_{i_{n+1}j_{n+1}}, \label{E:TwoElements}\\
\alpha_{2}= & \alpha_{i_{1}j_{1}}\alpha_{i_{2}j_{2}}\cdots
\alpha_{bc}\cdots\alpha_{(c+1)(d+1)}
\cdots \alpha_{i_{n+1}j_{n+1}}, \notag
\end{align}
giving $\alpha'$ via $(\partial^{c})^{*}$.

Let $a', a_i\in H_{2n}^{norm}(F(2n))$ be dual to the
elements
$\alpha', \alpha_i$.  We have thus shown that there are exactly four
elements
$\alpha \in H^{2n}_{norm}(F(2n))$ whose $d^{1}$, evaluated on $a'$ is
nonzero.  More precisely,
\begin{equation*}
\langle d^{1}\alpha_{1}, a'\rangle  =1, \ \
\langle d^{1}\alpha_{2}, a'\rangle  =-1, \ \
\langle d^{1}\alpha_{3}, a'\rangle  =-1, \ \
\langle d^{1}\alpha_{4}, a'\rangle  =1.
\end{equation*}
The signs depend on which index in
$a'$ is repeated (as $d^{1}$ is an alternating sum), as well as
indices $b, c,$ and $d$ (the number of transpositions needed to get
from one labeling to another has to be taken into account).  So
\begin{equation}\label{E:AlmostThere}
d_{1}a'=a_{1}-a_{2}-a_{3}+a_{4},
\end{equation}
and we conclude that
$$f\in ker(d^{1})\Longleftrightarrow \langle f,
a_{1}-a_{2}-a_{3}+a_{4}\rangle=0$$
for all $a_{1}$, $a_{2}$, $a_{3}$, and $a_{4}$ whose indices are
related as in \eqref{E:FourElements} or \eqref{E:TwoElements}.

If we recall that classes $a\in H_{2n}^{norm}(F(2n))$ were associated
to
chord diagrams with $n$ chords, then the conditions on $a_{i}$
precisely describe the four chord diagrams in the $4T$ relation.
Thus we get that, if $f$ is thought of as a function on
$CD_{n}$, then
$$
f\in E_{2}^{-2n, 2n}  \Longleftrightarrow \langle f, 4T\rangle
=0 \text{ for all $4T$ relations in $CD_{n}$}
                            \Longleftrightarrow f\in \W_{n}.
$$
\end{proof}

\subsection{The cotower of partial totalizations}
\label{S:Iso's&Surj'sInTower}

Recall that, in the normalized $E_{2}$ term of the spectral sequence,
we have
$$E_{2}^{-p,q}=0 \text{ \ if $p<q$ \ or \ $p=q=2k+1$ for some $k$.}$$
Since
$$E_{i+1}^{-p,q}=ker\left(E_{i}^{-p,q}\stackrel{d^{i}}
{\longrightarrow}E_{i}^{-p+i,q-i+1}\right)
\big/ im\left(E_{i}^{-p-i,q+i-1}\stackrel{d^{i}}{\longrightarrow}E_{i}^{-p,q}\right),$$
we have, for $p=q$,
\begin{equation}\label{E:DiagonalEntries}
E_{i+1}^{-p,p}=ker\left(E_{i}^{-p,p}\stackrel{d^{i}}
{\longrightarrow}E_{i}^{-p+i,p-i+1}\right).
\end{equation}

An immediate consequence is that we can construct an algebraic analog
of $T_{r}$, called $T_{r}^{*}$, as follows:  We truncate $\Xdot$
at the $r$th space for any $r$, and then construct the double
cochain complex (now finite in the horizontal direction) and the
associated spectral sequence.  Letting $T_{r}^{*}$ be the total
complex of this double complex, the spectral sequence would now be
computing the cohomology of $T_{r}^{*}$.
We then have
$$
H^{0}(T_{r}^{*})=\bigoplus_{j=0}^{r}E_{\infty}^{-j,j}.
$$
Now consider $H^{0}(T_{r+1}^{*})$.  A consequence of \eqref{E:DiagonalEntries}
is that its grading would be exactly the same as that of
$H^{0}(T_{r}^{*})$, but with one more summand.  In short, there would be
no new differentials coming into the diagonal of the spectral sequence
in passing from the truncation of $\Xdot$ at $r$ to $r+1$.  Since
every odd entry on the diagonal at the $E_{2}$ term is 0, we conclude
\begin{gather}
H^{0}(\Tstar{2n})\cong H^{0}(\Tstar{2n+1}) \label{E:TowerIso} \\
H^{0}(\Tstar{2n})/H^{0}(\Tstar{2n-1})\cong E_{\infty}^{-2n,2n}.
\label{E:TowerQuotient}
\end{gather}
We can put together the above isomorphisms (the second
giving an injection from
$
H^{0}(\Tstar{2n-1})$ to $H^{0}(\Tstar{2n})
$)
into a tower of invariants of partial totalizations which we call a
\emph{Taylor cotower}:

\begin{equation}\label{E:CohomologyTower}
\xymatrix{
\cdots \ar[r]^(0.3){\cong} &
H^{0}(\Tstar{2n-1}) \ar@{^{(}->}[r] &
H^{0}(\Tstar{2n}) \ar[r]^(0.45){\cong} &
H^{0}(\Tstar{2n+1}) \ar@{^{(}->}[r] &
\cdots
}
\end{equation}
Each space comes with a map $ev_{[i]}^{*}$ to $H^{0}(\overline{\K})$
described below.

\begin{rem}
Going from partial totalizations $T_{r}=Tot^{r}\Xdot$ to total
complexes $T_{r}^{*}=Tot^{r}C^{*}\Xdot$ is also possible
directly from the subcubical diagrams of configuration spaces.  We have
bypassed this by working with $\Xdot$ instead,  but we could have also
considered cochains on all spaces of punctured embeddings in the subcubical
diagram (all arrows would now get reversed)
and then formed the homotopy colimit $hocolim(C^{*}FC_{r})$ (dual to the homotopy limit, but
now in an algebraic sense).
A double complex $T_{r}^{*}$ could then be formed by collecting all
cochains on
embeddings with the same number of punctures.
\end{rem}

\subsection{The evaluation map to finite type invariants}\label{S:TotMapsToFinType}
Let $ev_{[i]}^{*}$ be
the map induced by $ev_{[i]}$ on cochains on configurations and cohomology:
$$
ev_{[i]}^{*}\colon
H^{0}(Tot^{i}C^{*}\Xdot)=H^{0}(T^{*}_{i})\longrightarrow
H^{0}(\overline{\K}).
$$
To simplify notation, we set
$
ev^{*}=ev_{[i]}^{*}
$
for all $i$.  Which value of $i$ is used will be clear from the
context throughout the section.

Now remember from previous section that
$
H^{0}(T^{*}_{2n})/H^{0}(T^{*}_{2n-1})=E^{-2n, 2n}_{\infty}.
$
\begin{prop}\label{P:TotMapsToFinType}
The image of
$
ev^{*}
$
is contained in $\V_{n}/\V_{n-1}$.  Further,
the diagram
$$
\xymatrix{ E^{-2n, 2n}_{\infty}
\ar[r]^{ev^{*}} \ar@{^{(}->}[d] &
\V_{n}/\V_{n-1} \ar@{^{(}->}[d]  \\
E^{-2n, 2n}_{2} \ar[r]^{\cong} & \W_{n}}
$$
commutes.
\end{prop}

For clarity, it is better to prove this proposition using the tower which is
dual to the one in \eqref{E:CohomologyTower}.  As mentioned in
\S\ref{S:Cd'sinSS}, we could have started with chains
on the cosimplicial space $\Xdot$ and formed the totalizations
$T^{i}_{*}=Tot^{i}C_{*}\Xdot$, where $Tot^{i}C_{*}\Xdot$ now means
the total complex of a double chain complex.  The tower of
totalizations then on homology looks like

\begin{equation}\label{E:HomologyTower}
\xymatrix{
\cdots &
H_{0}(T^{2n-1}_{*}) \ar[l]_(0.6){\cong} &
H_{0}(T^{2n}_{*}) \ar@{>>}[l] &
 H_{0}(T^{2n+1}_{*}) \ar[l]_(0.5){\cong} &
\cdots \ar@{>>}[l]
}
\end{equation}
There are also maps $ev_{*}$ from $H_{0}(\overline{\K})$ to each stage,
induced by $ev_{[i]}$ for each $i$ on
chains and homology.
Since we are ultimately interested in detecting finite type invariants in
the cotower, and,
consequently, in
mapping sums of resolutions of
singular knots into the stages of the homology tower, we
introduce an equivalent description of this tower.
\vskip 5pt
\noindent
Given a space $X$, one can consider the free abelian group $\R X$ of
linear combinations of points in $X$.  This is a space that inherits
its topology from $X$ and is related to the infinite symmetric
product of $X$.  The most important characterization of this space is
a version of the Dold-Thom Theorem, whose proof can be found in \cite{DT,H}:
\begin{thm}\label{T:Dold-Thom}  The homotopy groups of $\R X$ are
isomorphic to the homology groups of $X$ with real coefficients.
\end{thm}

The tower we will then
 use to prove \refP{TotMapsToFinType} is

\begin{equation}\label{E:HomotopyTower}
\xymatrix{
& \vdots \ar[d]^{\cong} \\
H_{0}(\overline{\K})=\pi_{0}(\R\overline{\K}) \ar[r]^{ev_{*}} \ar[dr]^{ev_{*}} \ar[ur] &
\pi_{0}(Tot^{2n}\R\Xdot) \ar@{>>}[d]   \\
&  \pi_{0}(Tot^{2n-1}\R\Xdot) \ar[d]^{\cong} \\
& \vdots}
\end{equation}
The advantage of this tower is that a singular knot can now be mapped into
it
by a map of the linear combination of its resolutions (with appropriate
signs as dictated by the Vassiliev skein relation).

\begin{lemma}
The towers \eqref{E:HomologyTower} and \eqref{E:HomotopyTower} are
equivalent.
\end{lemma}

\begin{proof}
Let
\begin{itemize}
\item $\mathcal{C}$ be the category of chain complexes,
\item $\mathcal{S}$ the category of spectra, and
\item $\mathcal{T}$ the category of topological spaces.
\end{itemize}
Given $C_{*}\in\mathcal{C}$, define $T^{i}C_{*}$ to be the truncation
of $C_{*}$,
$$
T^{i}C_{*}=ker(\partial)\stackrel{\partial}{\longleftarrow} C_{i+1}\stackrel{\partial}{\longleftarrow}
C_{i+2}\stackrel{\partial}{\longleftarrow}\cdots.
$$
A simplicial abelian group
$A^{i}$ can be associated to each $T^{i}C_{*}$ by the Dold-Kan
construction \cite{D,K}
(in fact, $\mathcal{C}$ and $\mathcal{A}$ are equivalent; see also
\cite{May}).  The realization functor $\vert\cdot\vert$ can then be applied to
$A^{i}$
to yield a space we denote by $|T^{i}C_{*}|$.  Further, one can then
pass to $\mathcal{S}$ by collecting these spaces as
$$
\big\{ |T^{0}C_{*}|,\  |\Sigma T^{-1}C_{*}|,\  |\Sigma^{2} T^{-2}C_{*}|,
\ldots\big\}.
$$
Here $\Sigma$ means an algebraic suspension which shifts the
truncation of $C_{*}$, so that the above spectrum is in fact an $\Omega$-spectrum.
We denote the resulting composed functor from $\mathcal{C}$ to $\mathcal{S}$
by $F$.

The homotopy of $|T^{0}C_{*}|$ can
be shown to be the same as the (non-negative) homology of $C_{*}$.
On the other hand, we can apply the functor
$\Omega^{\infty}$ to $F(C_{*})$ and get a space whose homotopy is
the homotopy of $|T^{0}C_{*}|$ since
$F(C_{*})$ is an $\Omega$-spectrum.

Now, given a space $X$, we can do the same:  Consider the composition
\begin{gather*}
\mathcal{T}\longrightarrow \mathcal{C} \longrightarrow
\mathcal{S}\longrightarrow \mathcal{T}
\\
X\longmapsto C_{*}X\longmapsto F(C_{*}X)\longmapsto
\Omega^{\infty}F(C_{*}X)\sim |T^{0}C_{*}X| \sim \R X.
\end{gather*}
The first equivalence in $\mathcal{T}$ is again due to the fact that $F(C_{*}X)$
is an $\Omega$-spectrum.  The second is given by \refT{Dold-Thom}
since the homotopy of $F(C_{*}X)$ is the homology of $|T^{0}C_{*}X|$,
which, in turn, is the homology of $X$.

But we can go even further, and start with a cosimplicial space $\Xdot$
(or any diagram of spaces viewed as a functor from an indexing
category to $\mathcal{T}$).  A generalized version of the above
produces a  cosimplicial
spectrum $F(C_{*}\Xdot)$ and we get a composition of functors
$$
\Xdot\longmapsto C_{*}\Xdot\longmapsto F(C_{*}\Xdot)\longmapsto
\Omega^{\infty}F(C_{*}\Xdot) \sim \R \Xdot.
$$
In particular, from the weak equivalence in above, we can deduce
\begin{equation}\label{E:TotF}
Tot^{r}(\Omega^{\infty}F(C_{*}\Xdot)) \sim Tot^{r}(\R
\Xdot),
\end{equation}
for all partial totalizations $Tot^{r}$.
\begin{rems}
Totalization only preserves equivalences, however, if $\Xdot$ is a
fibrant cosimplicial space.  Every cosimplicial space can be
replaced by a fibrant one via a degreewise equivalence of
cosimplicial spaces and in fact our cosimplicial space had to be
replaced by a fibrant one for \refT{Tot=holim} to be true anyway.

For a more general diagram of spaces, the role of the totalization is
played by the homotopy limit of the diagram.
\end{rems}
Next we have
\begin{equation}\label{E:TotLoopInfty}
\Omega^{\infty}F(Tot^{r}C_{*}\Xdot)
\sim
\Omega^{\infty}Tot^{r}(F(C_{*}\Xdot))
\sim
Tot^{r}(\Omega^{\infty}F(C_{*}\Xdot)).
\end{equation}
The second equivalence is true because the functor $Tot^{r}$
is defined as a compatible collection of spaces
$Maps(\Delta^{i}, X^{i})$, $i\leq r$, and
$
\Omega\, Maps(\Delta^{i}, X^{i})\sim Maps(\Delta^{i}, \Omega X^{i}).
$
Essentially the same reasoning applies in the case of the first
equivalence, except the argument is not as straightforward since
$Tot^{r}C_{*}\Xdot$ is an algebraic construction involving total
complexes rather than mapping spaces.

But we also know that there is an equivalence
\begin{equation}\label{E:TotC=TotR}
\Omega^{\infty}F(Tot^{r}C_{*}\Xdot)\sim Tot^{r}\R X.
\end{equation}
In particular, putting \eqref{E:TotF},
\eqref{E:TotLoopInfty}, and \eqref{E:TotC=TotR} together on $\pi_{0}$, we have
$$
\pi_{0}(Tot^{r}\R\Xdot)=H_{0}(Tot^{r}C_{*}\Xdot).
$$
Since the tower \eqref{E:HomologyTower} had
$H_{0}(Tot^{r}C_{*}\Xdot)=H_{0}(T_{*}^{r})$ as its stages, we can now
replace them by $\pi_{0}(Tot^{r}\R\Xdot)$.
\end{proof}

\begin{proof}[Proof of \refP{TotMapsToFinType}]
The idea will be that an $n$-singular knot maps to
$Tot^{2n-1}\R\Xdot$ via the evaluation map, but the map is homotopic
to 0.  However, the homotopy cannot be extended to $Tot^{2n}\R\Xdot$
and the obstruction to doing so is precisely the chord diagram
associated to the singular knot.

So let $K$ be an $n$-singular knot, and let $\{ a_{1}, \ldots , a_{n},
b_{1}, \ldots, b_{n}\}$ be the points in the interior of $I=[0, 1]$
making up the singularities, i.e.
$K(a_{i})=K(b_{i})$.  Let
$t=(t_{1}, \ldots, t_{j})$ be a point in $\Delta^{j}$ and note that
for all $j$ we are considering, $j<2n$.  The fact that there are
always fewer $t_{i}$ than $a_{i}$ and $b_{i}$ is crucial in what follows.

Now pick closed intervals $I_{k}$, $J_{k}$ around the
points $a_{k}$, $b_{k}$ so that the $2^{n}$ resolutions of $K$ differ
only in the $K(I_{k})$.  Let the resolutions be indexed by the subsets $S$
of $\{1, 2,\ldots , n\}$ in the sense that
if $k\in S$, then $K(I_{k})$ was changed into an overstrand (in the
Vassiliev skein relation).

It is clear how to extend $ev$ to linear combinations of knots
and get a map
$$
ev\colon \R\overline{\K} \ \longrightarrow \ (\Delta^{j} \to  \R X^{j}).
$$
In particular, we can map sums of resolutions $K_{S}$
(taken with appropriate signs) of an $n$-singular knot.

Let
\begin{gather*}
e_{k}^{+}=K_{S}(I_{k}), \ \ {\text{if}}\  k\in S, \\
 e_{k}^{-}=K_{S}(I_{k}), \ \ {\text{if}} \ k\notin S,
\end{gather*}
and define a homotopy (parametrized by $u$)
$$
e_{k}^{u}=(1-u)e_{k}^{+}+ue_{k}^{-}.
$$
Note that this homotopy is only defined on $I_{k}$ and has the effect of
changing an overcrossing to an undercrossing.  The problem is that,
for some time $u$, this homotopy reintroduces the original singularity on the
resolved knot.  This means that $ev$
will yield a degenerate configuration if any two of the $t_{i}$ equal
$a_{k}$ and $b_{k}$.

To avoid this, it suffices to modify the homotopy so that no movement takes
place when there is some $t_{i}\in\Delta^{j}$ which equals $b_{k}$.
To do this,
define a map
$$
\phi_{k}(t):\ \Delta^{j}\ \longrightarrow \ [0, 1]
$$ by
\[
\phi_{k}(t)=
\begin{cases}
\frac{2 min_{l}(|t_{l}-b_{k}|)}{|J_{k}|}, &  \text{ if at least one
of the $t_{i}$ is in $J_{k}$; } \\
 1 , & \text{ otherwise. }
\end{cases}
\]
Then $\phi_{k}(t)$ is a continuous map which is 0 if, for any
$i$, $t_{i}=b_{k}$ and 1 if, for all $i$, $t_{i}\notin J_{k}$.

Now let
\[
e_{S}^{u, t}=
\begin{cases}
e_{k}^{u\phi_{k}(t)}, &  \text{ in $I_{k}$ if $k\in S$; } \\
K_{S} , & \text{ outside $I_{k}$ for all $k$. }
\end{cases}
\]

Thus $e_{S}^{u, t}$ changes overcrossing $K_{S}(I_{k})$ to undercrossing
$K_{S}(J_{k})$ if
no $t_{i}$ is in $J_{k}$.  If some $t_{i}$
enters $J_{k}$, the homotopy happens only part of the way.
If any $t_{i}$ comes to the center of $J_{k}$ and thus equals
$b_{k}$, the homotopy is just constant.

This setup then prevents the collision of configuration points in
$X^{j}$ so that it makes sense to use $ev$ to define a homotopy
$$
K_{S}\times I\ \longrightarrow (\Delta^{j} \to X^{j})
$$
by
$$
(K_{S}, u)\ \mapsto \ ev(e_{S}^{u, t})(t).
$$
Adding over all $S\subset \{1, \ldots, n\}$ with appropriate signs, we get
a homotopy
\begin{equation}\label{E:homotopy}
\sum_{S}(\pm)(K_{S}\times I) \ \longrightarrow \ (\Delta^{j}\to \R X^{j}).
\end{equation}
We wish for this homotopy to be 0 when $u=1$.  It is important to keep in mind
that, for every $K_{S}$ and every $k$, there exists a $K_{S'}$ such
that $ev(K_{S})$ and $ev(K_{S'})$ differ only in those $t_{i}$ which
are in $I_{k}$, and these two configurations have opposite signs.

We pick a point $t\in\Delta^{j}$ and distinguish the following cases:

\vskip 5pt
\noindent
{\bf Case 1:}  $\forall i$ and $\forall k,$  $t_{i}\notin I_{k}$.

Then $ev(e_{S}^{u,t})(t)$ is the same configuration for all $S$ and
for all $u$, and
these cancel out in $\R X^{j}$.  Note that this happens regardless of whether
$J_{k}$ contains any of the $t_{i}$ or not.

\vskip 5pt
\noindent
{\bf Case 2:}  $\exists i$ and $\exists k$, such that $t_{i}\in I_{k}$.

{\bf Case 2a:}  $\forall j, \ t_{j}\notin J_{k}$.

Since $J_{k}$ is free of the $t_{j}$, $ev(e_{S}^{u, t})(t)$ moves the
overstrand $K_{S}(I_{k})$ to the understrand $K_{S}(I_{k})$ without
introducing a singularity.  So
now the two resolutions are the same and they cancel because of the
difference in signs.  Since all the resolutions can be combined in
pairs which only differ in one crossing, everything
adds up to 0.

{\bf Case 2b:}  $\exists j$ such that $t_{j}\in J_{k}$.

In this case, the homotopy $ev(e_{S}^{u, t})(t)$ does not change
the overcrossing to the undercrossing (it might perform a part of the movement,
depending on where exactly $t_{j}$ is in $J_{k}$---it is important
here that this partial movement happens on half the resolutions which
will pair up and cancel, as will the other half).  But since $j<2n$,
there must be an interval, either $I_{k'}$ or $J_{k'}$, which is free
of all the $t_{l}$.  This puts us back into one of the previous cases, with
$k=k'$.  Everything cancels in pairs again.

\vskip 5pt
\noindent
This exhausts all the cases so that
$$
\sum_{S}(\pm)K_{S} \ \longrightarrow \ (\Delta^{j} \to \R X^{j})
$$
is homotopic to 0 for $1\leq j \leq 2n-1$.  In other words, the map
$$
ev_{*}: \pi_{0}(\Z E) \ \longrightarrow \ \pi_{0}(\Delta^{j} \to \R X^{j})
$$
is 0 on resolutions of $n$-singular knots for
 $1\leq j \leq 2n-1$.

We still need this to be a homotopy in $Tot^{2n-1}\R\Xdot$,
so that we need to check the compatibility with coface and
codegeneracy maps in $\Xdot$.  This, however, is immediate from the
definition of $\phi_{k}(t)$ and we omit
the details.
Therefore
$$ev_{*}: \pi_{0}(\Z E) \ \longrightarrow \ \pi_{0}(Tot^{2n-1}\Z E)$$
is 0 on $n$-singular knots.

\vskip 5pt
\noindent
If we try to apply the same construction to $Tot^{2n}\R\Xdot$,
we will get 0 at $u=1$ if $j<2n$, but Case 2b now
breaks down when $j=2n$.  This is because the $t_{i}$ could land in {\it all}
$I_{k}$ and $J_{k}$.  Part of the homotopy might still take place
(depending on exactly where the $t_{i}$ are in those intervals), but
there will be no cancelations among the configurations.  We
distinguish an extreme case which will be needed later, namely let
$$
t^{*}=\{a_{1}, \ldots, a_{n}, b_{1}, \ldots, b_{n}\}\in\Delta^{2n}.
$$
However, since $ev(e_{S}^{u, t})(t)$ still produces 0 outside of
$$
B^{2n}=I_{1}\times J_{1}\times \cdots \times I_{n}\times J_{n} \in
int(\Delta^{2n}),
$$
we have, at $u=1$, a map
$
\Delta^{2n} \to  \R X^{2n}
$
which takes everything outside of
$B^{2n}$ to 0.  In particular, the boundary of $\Delta^{2n}$ goes to 0, so we
in fact have a map
$
S^{2n} \to \R X^{2n},
$
or an element of
$
\pi_{2n}(\R X^{2n})\cong H_{2n}(X^{2n}).
$
To see which element this is, we turn to cohomology:

Since $X^{2n}$ is the configuration space $F(2n, \R^{3})$,
we in fact have a map
$
B^{2n}\to \R F(2n, \R^{3})
$
and various maps
$
\R F(2n, \R^{3})\to \R S^{2}
$
given by the normalized difference of two points in the
configuration.  Since we are looking for an element of $H^{2n}$, we
take $n$ such maps to pull back a $2n$-form from a product of $n$
spheres.

Now remember that $B^{2n}$ is a product of $2n$ disjoint intervals,
so that the composition
$$
B^{2n}\to \R F(2n, \R^{3}) \to (\R S^{2})^{n}
$$
breaks up as $n$ maps
$$
I_{k}\times J_{k'} \to S^{2} \ \ \text{(or $I_{k}\times I_{k'}$ or
 $J_{k}\times J_{k'}$)}.
 $$
 But the boundary of each of these squares
 maps to 0, so we actually have $n$ maps from
 $S^{2}$ to $S^{2}.$
 Every such
 map is determined up to homotopy by its degree.  We will therefore first
pick a point $(v_{1}, \ldots, v_{n})$ in $(S^{2})^{n}$ which is a
 regular value for the composition
  $$
 B^{2n}\to \R F(2n, \R^{3}) \to (\R S^{2})^{n}.
 $$
 But since the linear combination of knots is here given by the
 resolutions of an $n$-singular knot, we may pick a point which is a
 regular value for
 $$
 B^{2n}\to F(2n, \R^{3}) \to (S^{2})^{n}
 $$
for all resolutions $S$.  Then we will count the number of preimages of this
regular value over all $S$ and all possible vectors between points in
the configuration lying in the images of intervals $I_{k}$
and $J_{k}$.  This will give an element of $H^{2n}(X^{2n})$ whose dual
is precisely the element which prevents the extension of the homotopy
to $Tot^{2n}\R\Xdot$.

Let $S=\{1, \ldots, n\}$ and choose $(v_{1}^{*}, \ldots, v_{n}^{*})\in (S^{2})^{n}$ such
 that $v_{i}^{*}$ is the (normalized) vector which points from $ev(K_{S})(b_{i})$
  to
 $ev(K_{S})(a_{i})$.  Then the only time the
 map $F(2n, \R^{3}) \to (S^{2})^{n}$ produces this value is when $t=t^{*}$.
 Also, this value occurs only once over all the resolutions, and that
 is when the singular knot is resolved positively at all the singularities (so
 $K(I_{k})$ is always the overstrand).
 Other simple geometric considerations show that
 this is the only preimage.

 We may thus think of the cohomology class we
 obtain as represented by $n$ vectors pointing from the images
 of the $a_{i}$ to the images of the $b_{i}$.  But then the dual class
 may be represented by a chord diagram which is
 precisely the diagram associated to the $n$-singular knot we started
 with.

Finally note that
$ev_{*}\colon \pi_{0}(\R \overline{\K}) \to \pi_{0}(Tot^{2n-1}\R \Xdot)$ is
0 for any singular knot with $n$ or more singularities, but it is not 0
on knots with $n-1$ singularities (repeat the arguments with a
shift---the key is that we now would not have $j<2(n-1)$ for all $j$).

\vskip 5pt
\noindent
We can now collect all the information and prove the proposition as
stated:
Let
\begin{align*}
\overline{\K}_{n} & =\{\text{sums of resolutions of $n$-singular knots}\} \\
[\R\overline{\K}_{n}] & =\{\text{homotopy classes of combinations of elements of
$\overline{\K}_{n}$}\}.
\end{align*}
Notice that there is a filtration
\begin{equation}\label{E:Filtration}
\pi_{0}(\R\overline{\K})\supset [\R\overline{\K}_{1}]\supset \cdots \supset [\R\overline{\K}_{n}] \supset
[\R\overline{\K}_{n+1}] \supset \cdots
\end{equation}
 since every
$(n+1)$-singular knot resolves into two $n$-singular ones.  Also
observe that
$$
\V_{n}=\big(\pi_{0}(\overline{\K})/[\R\overline{\K}_{n+1}]\big)^{*}
$$
and so
\begin{equation}\label{E:Duals}
\V_{n}/\V_{n-1}=\big([\R\overline{\K}_{n}]/[\R\overline{\K}_{n+1}]\big)^{*}.
\end{equation}
We have thus produced a map
$$
H_{0}(\overline{\K})=\pi_{0}(\R\overline{\K})\supset [\R\overline{\K}_{n}]
\longrightarrow ker\big(\pi_{0}(Tot^{2n}\R\Xdot)\to
\pi_{0}(Tot^{2n-1}\R\Xdot)\big).
$$
Dual to the statement from previous section that the quotient on cohomology is
$E_{\infty}^{-2n,2n}$, the above kernel is precisely
$E^{\infty}_{-2n,2n}$.  On the other hand, this is a subgroup of $E^{2}_{-2n,2n}$
which \refP{E2} shows (dually) to be isomorphic to $\cd{n}$.  Thus we
have a map
\begin{equation}\label{E:GettingToTheSquare}
[\R\overline{\K}_{n}]\longrightarrow E^{\infty}_{-2n,2n}\subseteq
E^{2}_{-2n,2n}\cong\cd{n}.
\end{equation}
The fact that our homotopy produces 0 in $\pi_{0}(Tot^{2n}\R\Xdot)$ on
singular knots with more than $n$ singularities means that
\eqref{E:GettingToTheSquare} actually maps from the combinations of
resolutions of $n$-singular knots which are not combinations of
resolutions of $(n+1)$-singular ones, so that we may write
$$
ev_{*}\colon [\R\overline{\K}_{n}]/[\R\overline{\K}_{n+1}]\longrightarrow E^{\infty}_{-2n,2n}.
$$
Using \eqref{E:Duals}, the dual of this map is precisely the one in
the statement of the Proposition:
$$
ev^{*}\colon E_{\infty}^{-2n,2n}\longrightarrow \V_{n}/\V_{n-1}.
$$
For the second part, recall that the chord diagram we obtained in
$E^{2}_{-2n,2n}$ was the same
one that produced the singular knot in the usual way.  It follows
immediately that we have
a commutative diagram
\begin{equation}\label{E:TheSquareOnHomology}
\xymatrix{
E^{\infty}_{-2n,2n}   &   [\R\overline{\K}_{n}]/[\R\overline{\K}_{n+1}] \ar[l]
 \\
E^{2}_{-2n,2n} \ar@{->>}[u] &  \cd{n} \ar@{->>}[u] \ar[l]_{\cong}
}
\end{equation}
The surjection on the right is simply the statement that every
$n$-singular knot can be obtained from a chord diagram by embedding an
interval while identifying the points paired off by chords.
The dual diagram is then
\begin{equation}\label{E:TheSquareOnCohomology}
\xymatrix{
E_{\infty}^{-2n,2n} \ar[r] \ar@{^{(}->}[d]  &  \V_{n}/\V_{n-1} \ar@{^{(}->}[d] \\
E_{2}^{-2n,2n} \ar[r]^{\cong} &  \W_{n}
}
\end{equation}
as desired.
\end{proof}
The first part of \refP{TotMapsToFinType} thus essentially states that
any knot invariant which factors through the tower has to be finite
type.  In particular, recall \refT{B-TMainTheorem} and
notice that $T(W)$ can be extended linearly to sums of resolutions
of singular knots, thereby fitting the composition
$$
\xymatrix{
[\R\overline{\K}_{n}] \ar[r]^(0.3){ev_{*}}  & \pi_{0}(Tot^{2n}\R
\Xdot)=H_{0}(T^{2n}_{*}) \ar[r]^(0.8){T(W)}  & \R.
}
$$
But we now know that this composition must vanish
on $[\R\overline{\K}_{n+1}]$, so that
we
immediately have
\begin{cor}\label{C:T(W,h)IsFinType}  The invariant $T(W)\colon
\pi_{0}(\HO_{2n})\longrightarrow \R$, when extended to
$$
\pi_{0}(holim(\R
EC_{2n}))=\pi_{0}(Tot^{2n}\R \Xdot)=H_{0}(T^{2n}_{*}),
$$
vanishes on points in $T^{2n}_{*}$ which come from sums of
resolutions of $(n+1)$-singular
knots.  Thus $T(W)$ is a type $n$ invariant.
\end{cor}

\subsection{Proofs of the main theorems}\label{S:Universality}

We finally prove Theorems
\ref{T:IntroMainTheorem1} and \ref{T:Collapse}, thus in
particular
showing
that the algebraic analog
\eqref{E:CohomologyTower} of the Taylor
tower for $\overline{\K}$ classifies finite type invariants.

\vskip 5pt
\noindent
Since the Bott-Taubes integral
$ T(W)\in
H^{0}(T^{*}_{2n}), $
applied to a point $h(K)\in\HO_{2n}$ coming
from an $n$-singular knot $K$, is 0 in $H^{0}(T^{*}_{2n-1})$, we
think of it as a map
$$
T(W)\colon \W_{n}\longrightarrow
H^{0}(T^{*}_{2n})/H^{0}(T^{*}_{2n-1}),
$$
or, in other words, a map $\W_{n}\to E_{\infty}^{-2n,2n}$ on $n$-singular knots.
It thus fits the diagram
\eqref{E:TheSquareOnCohomology} as
\begin{equation}\label{E:CohomologySquareWithT}
\xymatrix{
E_{\infty}^{-2n,2n} \ar@{^{(}->}[rr] \ar@{^{(}->}[d]  &  &
\V_{n}/\V_{n-1} \ar@{^{(}->}[d] \\
E_{2}^{-2n,2n} \ar[rr]^{\cong}  & & \W_{n} \ar[ull]_(0.4){T(W)}
}
\end{equation}
We have also modified the diagram by noticing that the top map has to
be an injection
because of the commutativity of the
diagram.

But now we can combine \refT{UniversalInvariant} and the fact from
\refT{B-TMainTheorem} that $W(D)$ restricts to a universal finite
type invariant of ordinary knots.  In particular, the following is easy to see:

\begin{prop}\label{P:UniversalInvariant}
The composition
\begin{equation}\label{E:CompositionIsIdentity}
\xymatrix{
\W_{n} \ar[r] & E_{\infty}^{-2n,2n} \ar@{^{(}->}[r] & \V_{n}/\V_{n-1}
\ar@{^{(}->}[r] & \W_{n}
}
\end{equation}
is the identity.
\end{prop}
The first observation now is that the map
$$
\xymatrix{ \V_{n}/\V_{n-1} \ar@{^{(}->}[r]  & \W_{n}}
$$
must be a surjection, so that we have obtained an alternative proof
of the
Kontsevich Theorem:
\begin{thm}\label{T:NewKontsevich}$\V_{n}/\V_{n-1} \cong \W_{n}$.
\end{thm}

We then also have
\begin{equation}\label{E:QuotientIso}
H^{0}(T^{*}_{2n})/H^{0}(T^{*}_{2n-1})=E_{\infty}^{-2n,2n}
\cong \V_{n}/\V_{n-1}.
\end{equation}
But now we can compare two exact sequences
$$
\xymatrix{
0 \ar[r] & H^{0}(T^{*}_{2n-1}) \ar[r] \ar[d] &  H^{0}(T^{*}_{2n})
\ar[r] \ar[d] &
H^{0}(T^{*}_{2n})/H^{0}(T^{*}_{2n-1})  \ar[r] \ar[d]_{\cong} &  0 \\
0 \ar[r] & \V_{n-1} \ar[r] & \V_{n} \ar[r] & \V_{n}/\V_{n-1} \ar[r] &
0
}
$$
The first two vertical maps are duals to the ones given in
\eqref{E:GettingToTheSquare}.

Since the right square commutes, \eqref{E:QuotientIso} can be used
with $n=1$ as the base case for an induction which shows that the
middle vertical map is an isomorphism if the left one is.  We thus
have

\begin{thm}\label{T:MainTheorem1} $H^{0}(T^{*}_{2n})\cong \V_{n}.$
\end{thm}
\noindent
Finally,
$\xymatrix{E_{\infty}^{-2n,2n}\ar@{^{(}->}[r] & E_{2}^{-2n,2n} }$
also has
to be an isomorphism since it is the last map in the
commutative square with three other isomorphisms.  We rephrase
this as
\begin{thm}\label{T:MainTheorem2}  The spectral sequence for the
cohomology of $\Xdot$ collapses at $E_{2}$ on the diagonal.
\end{thm}
\noindent
The preceding two results are precisely Theorems
\ref{T:IntroMainTheorem1}
and \ref{T:Collapse}.

\subsection{Some further questions}\label{S:Future}

The Taylor tower for $\overline{\K}$ is thus a potentially rich
source of information about finite type theory.  Its advantage is
that the stages $\HO_{r}$ lend themselves to purely topological
examination.  At the time, not much is known about these spaces,
but our results give evidence that a closer look may shed some
light on classical knot theory.  In particular:

\begin{itemize}

\item For knots in $\R^{n}$, $n>3$, the cohomology spectral
sequence is guaranteed to converge to the totalization of
$X^{\bullet}$. Combining this with the Goodwillie-Klein-Weiss
result stating that in this case the Taylor tower converges to the
space of knots, one can extract useful information about the
homotopy and cohomology of spaces of knots \cite{Dev, Pasc}.
However, this may not be true in the case $n=3$ because it is not
immediately clear that $C^{*}$ commutes with totalization. If the
spectral sequence indeed does not converge to the desired
associated graded, then potentially an even more interesting
question arises: What can we say about invariants pulled back from
the Taylor tower without passing to cochains first?  It is quite
possible that one in this case also obtains only finite type
invariants, but the genuine Taylor tower for $\overline{\K}$ could
also contain more information.

\item Can one gain new insight into finite type theory by studying the
homotopy types of spaces $\HO_{r}$?  In particular, can one gain
topological insight into the common thread between knots in $\R^{3}$
and $\R^{n}$, namely the Kontsevich and Bott-Taubes integral
constructions of finite type invariants?  This is to be expected,
since the latter type of integrals plays a
crucial role in our proofs.

\item What is the relationship between the Taylor tower and
Vassiliev's approach
to studying $\K$ through
simplicial spaces \cite{Vas}?

\item Do other finite type invariants (of homology spheres, for
example) also factor through constructions coming from calculus of
functors?  One should at least be able to carry the ideas
and constructions here from knots to braids without too much
difficulty.

\item Can one say anything useful about the inverse limit of the
Taylor tower?  The question of invariants of this limit surjecting
onto knot invariants is precisely the question of finite type
invariants separating knots.

\item   Recall that the
evaluation map $ev^{*}$ provides the isomorphism of
\refT{MainTheorem1}.  Budney, Conant, Scannell, and Sinha \cite{BCSS}
examine the
evaluation map closely related to ours to give a geometric
interpretation of the unique (up to framing) type 2 invariant.
Can one in general understand the geometry of finite type
invariants using the evaluation map?

\item As predicted by a result in \cite{GK}, is the map
 $\K\! \to\! \HO_{r}$ 0-connected, i.e. surjective on
 $\pi_{0}$?

\end{itemize}

\end{document}